\documentclass[12pt,leqno]{article}
\usepackage{amsfonts}
\usepackage{amsmath, amsthm, amsfonts, amssymb, color}
\usepackage{mathrsfs}
\renewcommand{\bar}{\overline}
\renewcommand{\hat}{\widehat}
\renewcommand{\tilde}{\widetilde}

\theoremstyle{definition}

\definecolor{wco}{rgb}{0.5,0.2,0.3}

\numberwithin{equation}{section} \theoremstyle{remark}

\newcommand{\ua}{\uparrow}
\newcommand{\da}{\downarrow}

\title{{\bf Distribution dependent BSDEs driven by Gaussian processes}
}
\author{
{\bf  Xiliang Fan$^{a)}$,  Jiang-Lun Wu$^{b)}$
 }\\
\footnotesize{$^{a)}$School of Mathematics and Statistics, Anhui Normal University, Wuhu 241002, China}\\
\footnotesize{$^{b)}$Department of Mathematics, Swansea University, Bay Campus, SA1 8EN, UK}\\
\footnotesize{\sf{fanxiliang0515@163.com},\ \sf{j.l.wu@swansea.ac.uk}}\\
}

\begin{document}
\def\R{\mathbb R}
\def\N{\mathbb N}
\def\E{\mathbb E}
 \def\H{\mathbb H}
\def\Q{\mathbb Q}
\def\H{\mathbb{H}}
\def\P{\mathbb{P}}
\def\S{\mathbb{S}}
\def\Y{\mathbb{Y}}
\def\W{\mathbb{W}}
\def\D{\mathbb{D}}

\def\cH{\mathcal{H}}
\def\cS{\mathcal{S}}

\def\sA{\mathscr{A}}
 \def\sB{\mathscr {B}}
 \def\sC{\mathscr {C}}
  \def\sD{\mathscr {D}}
 \def\sF{\mathscr{F}}
\def\sG{\mathscr{G}}
\def\sL{\mathscr{L}}
\def\sP{\mathscr{P}}
\def\sS{\mathscr{S}}
\def\sM{\mathscr{M}}
\def\eq{\equation}
\def\beg{\begin}
\def\ep{\epsilon}
\def\ve{\varepsilon}
\def\vp{\varphi}
\def\vr{\varrho}
\def\om{\omega}
\def\Om{\Omega}
\def\si{\sigma}
\def\ff{\frac}
\def\sq{\sqrt}
\def\kk{\kappa}
\def\de{\delta}
\def\<{\langle} 
\def\>{\rangle}
\def\Ga{\Gamma}
\def\ga{\gamma}
\def\na{\nabla}
\def\be{\beta}
\def\al{\alpha}
\def\pp{\partial}
 \def\ti{\tilde}
\def\1{\lesssim}
\def\ra{\rightarrow}
\def\da{\downarrow}
\def\upa{\uparrow}
\def\l{\ell}
\def\8{\infty}
\def\3{\triangle}
 \def\DD{\Delta}
\def\m{{\bf m}}
\def\B{\mathbf B}
\def\e{\text{\rm{e}}}
\def\la{\lambda}
\def\th{\theta}

\def\d{\text{\rm{d}}}
\def\ess{\text{\rm{ess}}}
\def\Ric{\text{\rm{Ric}}}
\def \Hess{\text{\rm{Hess}}}
\def\ua{\underline a}
\def\Ric{\text{\rm{Ric}}}
\def\cut{\text{\rm{cut}}}      
\def\alphaa{\mathbf{r}}     
\def\r{r}
\def\gap{\text{\rm{gap}}} 
\def\prr{\pi_{{\bf m},\varrho}}  
\def\r{\mathbf r}
\def\Tilde{\tilde} 
\def\TILDE{\tilde}
\def\II{\mathbb I}
\def\i{{\rm in}}
\def\Sect{{\rm Sect}}

\renewcommand{\bar}{\overline}
\renewcommand{\hat}{\widehat}
\renewcommand{\tilde}{\widetilde}

\allowdisplaybreaks

\maketitle
\begin{abstract}
In this paper we are concerned with distribution dependent backward stochastic differential equations (DDBSDEs) driven
by Gaussian processes. We first show the existence and uniqueness of solutions to this type of equations. This is done by
formulating a transfer principle to transfer the well-posedness problem to an auxiliary DDBSDE driven by Brownian motion.
Then, we establish a comparison theorem under Lipschitz condition and boundedness of Lions derivative imposed on the
generator. Furthermore, we get a new representation for DDBSDEs driven by Gaussian processes, this representation is
even new for the case of the equations driven by Brownian motion. The new obtained representation enables us to prove
a converse comparison theorem. Finally, we derive transportation inequalities and Logarithmic-Sobolev inequalities via the
stability of the Wasserstein distance and the relative entropy of measures under the homeomorphism condition.
\end{abstract}


\section{Introduction}

Backward stochastic differential equations (BSDEs) were first introduced in their linear form by Bismut in \cite{Bismut73} to investigate stochastic
control problems and their connections with a stochastic version of the Pontryagin maximum principle.
Afterwards, BSDEs was generally formalised and developed in the seminal work \cite{PP90}.
In the last decades, BSDEs have been the subject of growing interest in stochastic analysis,
as these equations naturally arise in stochastic control problems in mathematical finance and they provide Feynman-Kac type formulas for
semi-linear PDEs (see, e.g., \cite{EPQ9705,MY99,PP92,ZJ17}).

On the other hand, distribution dependent stochastic differential equations (DDSDEs), also known as mean-field equations or McKean-Vlasov
equations, are the It\^{o} equations whose coefficients depend upon the law of the solution. As DDSDEs can provide a probabilistic representation
for the solutions of a class of nonlinear PDEs, in which a typical example is the propagation of chaos, they are widely used as models in statistical
physics and in the study of large scale social interactions within the memory of mean-field games, for which we refer to, e.g., \cite{BT97,CD15,HMC13}
and references therein. Furthermore, nonlinear DDBSDEs were first introduced by Buckdahn, Djehiche, Li and Peng in \cite{BDLP09}.
Since then, the DDBSDEs have received increasing attentions and have been investigated in a variety of settings. Let us just mention a few here.
Chassagneux, Crisan and Delarue \cite{CCD15} showed the existence and uniqueness of solutions to fully coupled DDBSDEs; Carmona and Delarue \cite{CD15} studied DDBSDEs via the stochastic maximum principle; Li \cite{Li18} considered the well-posedness problem of DDBSDEs driven by a
Brownian motion and an independent Poisson random measure, and provided a probabilistic representation for a class of nonlocal PDEs of mean-field
type; Li, Liang and Zhang \cite{LLZ18} obtained a comparison theorem for DDBSDEs.

In this paper, we want to study the following DDBSDEs driven by Gaussian processes
\begin{equation}\label{Bsde-In}
\left\{
\begin{array}{ll}
\d Y_t=-f(t,X_t,Y_t,Z_t,\sL_{(X_t,Y_t,Z_t)})\d V_t+Z_t\d^\diamond X_t,\\
Y_T=g(X_T,\sL_{X_T}),
\end{array} \right.
\end{equation}
where $X$ is a centered one-dimensional Gaussian process such that $V_t:=\mathrm{Var}X_t, t\in[0,T]$, is a strictly increasing, continuous function with $V(0)=0$ introduced in \cite{Bender14}, $\sL_{(X_t,Y_t,Z_t)}$ and $\sL_{X_T}$ denote respectively the laws of $(X_t,Y_t,Z_t)$ and $X_T$,
and the stochastic integral is the Wick-It\^{o} integral defined by the $S$-transformation and the Wick product (see Section 2.1).
Precise assumptions on the generator $f:[0,T]\times\R\times\R\times\R\times\sP_\th(\R\times\R\times\R)\ra\R$ and the terminal value function $g:\R\times\sP_\th(\R)\ra\R$ will be specified in later sections, where $\sP_\th(\R^m)$ stands for the totality of probability measures on $\R^m$ with finite $\th$-th moment. We would like to mention that the driving noise $X$ in \eqref{Bsde-In} includes fractional Brownian motion $B^H$ with Hurst parameter $H\in(0,1)$ (while $B^{1/2}$ is the standard Brownian motion), and fractional Wiener integral (see Remark \ref{Re(Inte)}).

The main objectives of the present paper are to show the well-posedness and (converse) comparison theorems, and then to establish functional inequalities including transportation inequalities and Logarithmic-Sobolev inequalities for \eqref{Bsde-In}. Our strategy is as follows.
Based on a new transfer principle that extends \cite[Theroem 3.1]{Bender14} to the distribution dependent setting,
we first prove general existence and uniqueness results for \eqref{Bsde-In} (see Theorem \ref{Th1}), which are then applied to the case of Lipschitz
generator $f$.
Second, with the help of a formula for the $L$-derivative, we are able to derive a comparison theorem
which generalises and improves the corresponding one in the existing literature (see Theorem \ref{Th(com)} and Remark \ref{Re-comp1}).
It is worth stressing that comparing with the works in the distribution-free cases, here we need to impose an additional condition on $f$ that involves
the Lions derivative of $f$.
Moreover, we obtain a converse comparison theorem which is roughly speaking a converse to the comparison theorem obtained above (see Theorem \ref{Th(Conve)} and Remark \ref{Re(Conve)}).
To this end, we provide a representation theorem for the generator $f$ which is even new for DDBSDEs driven by Brownian motions which is also
interesting in itself. Finally, by utilising the stability of the Wasserstein distance and relative entropy of measures under the homeomorphism,
we establish several functional inequalities including transportation inequalities and Logarithmic-Sobolev inequalities (see Theorems \ref{Th(TrIn)} and \ref{Th(LS)}). Here, let us point out that the transportation inequality for the law of the control solution $Z$ of \eqref{Bsde-In} stated in
Theorem \ref{Th(TrIn)} is of the form
\begin{align*}
\W_p(\sL_{Z},\mu)\leq C\left(H(\mu|\sL_{Z})\right)^{\ff 1 {2p}}
\end{align*}
with any $p\geq1$.
In particular, when $p=2$, this inequality reduces to
\begin{align*}
\W_2(\sL_{Z},\mu)\leq C\left(H(\mu|\sL_{Z})\right)^{\ff 1 {4}},
\end{align*}
which is clearly different from the usual quadratic transportation inequality (also called the Talagrand inequality). On the other hand,
as shown in \cite{BT20}, this type of inequalities allows one to derive deviation inequality. However, it does not allow to get other important inequalities
such as Poincar\'{e} inequality etc. Hence, an interesting problem is whether our results can be further improved in the sense of $\W_2(\sL_{Z},\mu)\leq C\sqrt{H(\mu|\sL_{Z})}$. Our techniques are currently not enough to give a full answer, since Lemma \ref{FI-Le2} below cannot applied to the case of $Z$.
We will leave this topic for the future work.

The remaining of the paper is organised as follows. Section 2 presents some basic facts on Gaussian processes, the Lions derivative,
and introduce a transfer result which allows us to build a relation between DDBSDEs concerned and DDBSDEs driven by Brownian motion.
In Section 3, we show the existence and uniqueness of a solution to DDBSDE driven by Gaussian process.
In Section 4, we establish a comparison theorem, and also provide a converse comparison theorem via a new representation theorem.
Section 5 is devoted to deriving functional inequalities, including transportation inequalities and Logarithmic-Sobolev inequalities.
Section 6 is designed as an appendix that we prove an auxiliary result needed in Section 5 (cf. Proposition 5.1).

\section{Preliminaries}

\subsection{Wick-It\^{o} integral for Gaussian processes}

In this part, we shall recall some important definitions and facts concerning the Wick-It\^{o} integral for Gaussian processes.
Further detailed and deep discussions can be found, e.g., \cite[Section 2]{Bender14} and references therein.

Let $(\Omega,\sF,(\sF_t^X)_{t\in[0,T]},\P)$ be a filtered probability space with $(\sF_t^X)_{t\in[0,T]}$ the natural completed and right continuous filtration generated by a centered Gaussian process $(X_t)_{0\leq t\leq T}$,
whose covariance function $V_t:=\mathrm{Var}X_t, t\in[0,T]$, is a strictly increasing and continuous function with $V(0)=0$.

The first chaos associated to $X$ is
\begin{align*}
\mathfrak{C}_X:=\overline{\mathrm{span}\{X_t:t\in[0,T]\}},
\end{align*}
where the closure is taken in $L^2_X:=L^2(\Omega,\sF_T^X,\P)$.
It is obvious that the elements in $\mathfrak{C}_X$ are centered Gaussian variables.
We define the map $\mathcal{R}:\mathfrak{C}_X\ra\R^{[0,T]}$ by
\begin{align*}
(\mathcal{R}f)(t)=\E(X_tf).
\end{align*}
It is readily checked that $\mathcal{R}$ is injective, whose image $\mathcal{R}(\mathfrak{C}_X)$ is called the Cameron-Martin space of $X$ and is equipped with the inner product
\begin{align*}
\langle f,g\rangle_X:=\E\left[\mathcal{R}^{-1}(f)\mathcal{R}^{-1}(g)\right].
\end{align*}
Now, we let $\mathfrak{H}_X$ be the set of all those $\mathfrak{h}\in\mathcal{R}(\mathfrak{C}_X)$ that are absolutely continuous with respect to  $\d V$ with square integrable density, i.e.
\begin{align*}
\mathfrak{h}(t)=\int_0^t\dot{\mathfrak{h}}(s)\d V_s, \ \ \dot{\mathfrak{h}}\in L^2([0,T],\d V).
\end{align*}
Throughout the paper, we suppose that $\mathfrak{H}_X$ and $\{\dot{\mathfrak{h}}:\mathfrak{h}\in\mathfrak{H}_X\}$ are respectively the dense subsets of
$\mathcal{R}(\mathfrak{C}_X)$ and $L^2([0,T],\d V)$.

\beg{rem}\label{Re(Inte)}
As pointed out in \cite[Theorem 2.2]{Bender14}, the Gaussian processes concerned consists of a large class of examples which includes, e.g.,
fractional Brownian motion $B^H$ with $H\in(0,1)$ and fractional Wiener integral of the form $\int_0^t\si(s)\d B^H_s$ with $H\in(1/2,1)$
and a deterministic function $\si$ satisfying $c^{-1}\leq\si\leq c$ for some $c>0$.
\end{rem}

Next, we shall introduce the construction of the Wick-It\^{o} integral with respect to Gaussian process $X$.
Due to \cite[Corollary 3.40]{Janson97}, the random variables of the form
\begin{align*}
\e^{\diamond\mathfrak{h}}:=\exp\left\{\mathcal{R}^{-1}(\mathfrak{h})-\ff 1 2 \mathrm{Var}\mathcal{R}^{-1}(\mathfrak{h})\right\},\ \ \mathfrak{h}\in\mathfrak{H}_X,
\end{align*}
form a total subset of $L^2_X$, in which $\e^{\diamond\mathfrak{h}}$ is called Wick exponential.
Then for each random variable $\eta\in L^2_X$, it can be uniquely determined by its $\mathcal{S}$-transform
\begin{align*}
(\mathcal{S}\eta)(\mathfrak{h}):=\E(\eta\e^{\diamond\mathfrak{h}}), \ \ \mathfrak{h}\in\mathfrak{H}_X.
\end{align*}
That is, if $\eta$ and $\zeta$ belong to $L^2_X$ satisfying $(\mathcal{S}\eta)(\mathfrak{h})=(\mathcal{S}\zeta)(\mathfrak{h})$ for every $\mathfrak{h}\in\mathfrak{H}_X$, then there holds $\eta=\zeta, \P$-a.s..
In addition, observe that for each $\mathfrak{h}\in\mathfrak{H}_X$,
\begin{align*}
(\mathcal{S}X_t)(\mathfrak{h})=\E(X_t\mathcal{R}^{-1}(\mathfrak{h}))=\mathfrak{h}(t)=\int_0^t\dot{\mathfrak{h}}(s)\d V_s, \ \ t\in[0,T]
\end{align*}
is a bounded variation function and then can be regarded as an integrator in a Lebesgue-Stieltjes integral, which allows us to introduce the following  Wick-It\^{o} integral.

\beg{defn}\label{De-WI}
A measurable map $Z:[0,T]\ra L^2_X$ is said to have a Wick-It\^{o} integral with respect to $X$, if for any $\mathfrak{h}\in\mathfrak{H}_X$,
\begin{align*}
\int_0^T(\mathcal{S}Z_t)(\mathfrak{h})\d \mathfrak{h}(t)
\end{align*}
exists and there is a random variable $\xi\in L^2_X$ such that
\begin{align*}
(\mathcal{S}\xi)(\mathfrak{h})=\int_0^T(\mathcal{S}Z_t)(\mathfrak{h})\d \mathfrak{h}(t).
\end{align*}
In this case, we denote $\xi$ by $\int_0^TZ_t\d^\diamond X_t$ and call it the Wick-It\^{o} integral of $Z$ with respect to $X$.
Besides, we often use $\int_a^bZ_t\d^\diamond X_t$ to denote $\int_0^T\mathrm{I}_{[a,b]}(t)Z_t\d^\diamond X_t$.
\end{defn}

\beg{rem}\label{Re(WI)}
(i) Suppose that $Z:[0,T]\ra L^2_X$ is continuous and $\{\pi^n\}_{n\geq1}$ is a sequence of partition of $[0,T]$. Then, we have
\begin{align}\label{1Re(WI)}
\int_0^TZ_t\d^\diamond X_t=\lim_{|\pi^n|\ra0}\sum_{t_i\in\pi^n}Z_{t_i}\diamond(X_{t_{i+1}}-X_{t_i}),
\end{align}
provided that the above limit exists in $L^2_X$.
Here, $Z_{t_i}\diamond(X_{t_{i+1}}-X_{t_i})$ is a Wick product defined as follows:
\begin{align*}
(\mathcal{S}(Z_{t_i}\diamond(X_{t_{i+1}}-X_{t_i})))(\mathfrak{h})=(\mathcal{S}Z_{t_i})(\mathfrak{h})(\mathcal{S}(X_{t_{i+1}}-X_{t_i}))(\mathfrak{h}),\ \ \mathfrak{h}\in\mathfrak{H}_X.
\end{align*}
In view of \eqref{1Re(WI)}, one can see that the Wick-It\^{o} integral can be interpreted as a limit of Riemann sums in terms of the Wick product.

(ii) If $X=B^{1/2}$, i.e. $X$ is a Brownian motion, then we obtain that $V_t=t$ and
\begin{align*}
\mathcal{R}(\mathfrak{C}_{B^{1/2}})=\mathfrak{H}_{B^{1/2}}=\left\{\mathfrak{h}:\mathfrak{h}(t)=\int_0^t\dot{\mathfrak{h}}(s)\d s, \ \ \dot{\mathfrak{h}}\in L^2([0,T],\d t)\right\}.
\end{align*}
Suppose that $Z$ is progressively measurable satisfying $\E\int_0^TZ_t^2\d t<\infty$. Then, it is easy to verify that the Wick-It\^{o} integral $\int_0^TZ_t\d^\diamond B^{1/2}_t$ coincides with the usual It\^{o} integral $\int_0^TZ_t\d B^{1/2}_t$, and then
\begin{align}\label{2Re(WI)}
\left(\mathcal{S}\int_0^TZ_t\d B^{1/2}_t\right)(\mathfrak{h})=\int_0^T(\mathcal{S}Z_t)(\mathfrak{h})\d \mathfrak{h}(t),\ \ \mathfrak{h}\in\mathcal{R}(\mathfrak{C}_{B^{1/2}})=\mathfrak{H}_{B^{1/2}}.
\end{align}

More details can be found in \cite[Remark 2.4]{Bender14}.
\end{rem}

\subsection{The transfer principle}

This part is devoted to establishing a transfer principle, which connects DDSDEs driven by Gaussian processes and DDSDEs driven by Brownian motion,
and will play a crucial role in the proofs of our main results.

In order to state the principle, we let $U$ be the inverse of $V$ defined as
\beg{align*}
U_s:=\inf\{r\geq0: V_r\geq s\}, \ \ s\in[0,V_T],
\end{align*}
and introduce an auxiliary Brownian motion $(\ti{W})_{t\in[0,V_T]}$ on a filtered probability space $(\ti\Omega,\ti\sF,(\ti\sF_t^X)_{t\in[0,V_T]},\ti\P)$,  is the filtration generated by $\ti{W}$.
Similar to Section 2.1, we can define the $\mathcal{S}$-transform on this auxiliary probability space as follows:
for each random variable $\ti\eta\in L^2(\ti\Omega,\ti\sF,\ti\P)$,
\begin{align*}
(\tilde{\mathcal{S}}\ti\eta)(\mathfrak{h}):=\ti\E\left(\ti\eta\exp\left\{\int_0^{V_T}\dot{\mathfrak{h}}(s)\d\ti{W}_s-\ff 1 2\int_0^{V_T}\dot{\mathfrak{h}}^2(s)\d s\right\}\right), \ \ \mathfrak{h}\in\mathcal{R}(\mathfrak{C}_{\ti{W}}).
\end{align*}
Here, we recall that owing to Remark \ref{Re(WI)} (ii), $\mathcal{R}(\mathfrak{C}_{\ti{W}})$ is of the form
\begin{align*}
\mathcal{R}(\mathfrak{C}_{\ti{W}})=\mathfrak{H}_{\ti{W}}=\left\{\mathfrak{h}:\mathfrak{h}(t)=\int_0^t\dot{\mathfrak{h}}(s)\d s, \ \ \dot{\mathfrak{h}}\in L^2([0,V_T],\d t)\right\}.
\end{align*}

Now, we have the following transfer principle, which is a distribution dependent version of \cite[Theorem 3.1]{Bender14}.

\beg{prp}\label{Pr1}
Assume that $\phi:\R\times\sP_2(\R)\ra\R, b:[0,T]\times\R\times\sP_2(\R^m)\ra\R, \vartheta:[0,T]\times\R\times\sP_2(\R)\ra\R^m$ and $\si:[0,T]\times\R\times\sP_2(\R)\ra\R$
are measurable functions satisfying
\beg{align}\label{Pr1-0}
\E\left[\phi^2(X_t,\sL_{X_t})+\int_0^t(b^2(s,X_s,\sL_{\vartheta(s,X_s,\sL_{X_s})})+\si^2(s,X_s,\sL_{X_s}))\d V_s\right]<\infty \end{align}
for some $t\in[0,T]$, and
\beg{align}\label{Pr1-1}
\phi(\ti{W}_{V_t},\sL_{\ti{W}_{V_t}})=\int_0^{V_t}b(U_s,\ti{W}_s,\sL_{\vartheta(U_s,\ti{W}_s,\sL_{\ti{W}_s})})\d s+\int_0^{V_t}\si(U_s,\ti{W}_s,\sL_{\ti{W}_s})\d\ti{W}_s, \ \ \tilde{\P}\textit{-}a.s.
\end{align}
Then, $\int_0^t\si(s,X_s,\sL_{X_s})\d^\diamond X_s$ is well-defined and there holds in $L^2(\Omega,\sF_T^X,\P)$
\beg{align}\label{Pr1-2}
\phi(X_t,\sL_{X_t})=\int_0^tb(s,X_s,\sL_{\vartheta(s,X_s,\sL_{X_s})})\d V_s+\int_0^t\si(s,X_s,\sL_{X_s})\d^\diamond X_s.
\end{align}
\end{prp}

Before proving Proposition \ref{Pr1}, we first give a useful lemma whose proof is identical to \cite[Lemma 3.2]{Bender14} and therefore omitted here.

\beg{lem}\label{Le1}
Assume that $\vartheta:[0,T]\times\R\times\sP_2(\R)\ra\R^m$ and $\psi:[0,T]\times\R\times\sP_2(\R^m)\ra\R$ are two measurable functions
such that $\E\psi^2(t,X_t,\sL_{\vartheta(t,X_t,\sL_{X_t})})<\infty$ with some $t\in[0,T]$.
Then for any $\hbar\in\mathfrak{H}_X$,
\beg{align}\label{1Le1}
(\mathcal{S}\psi(t,X_t,\sL_{\vartheta(t,X_t,\sL_{X_t})}))(\hbar)=(\tilde{\mathcal{S}}\psi(t,\tilde{W}_{V_t},\sL_{\vartheta(t,\tilde{W}_{V_t},\sL_{\tilde{W}_{V_t}})}))(\hbar\circ U).
\end{align}
\end{lem}
Let us stress that if $\hbar\in\mathfrak{H}_X$, then $\hbar\circ U$ belongs to $\mathcal{R}(\mathfrak{C}_{\ti{W}})$, and thus the right-hand side of \eqref{1Le1} is well-defined.
Indeed, observe that
\beg{align*}
(\hbar\circ U)(t)=\int_0^{U_t}\dot{\hbar}(s)\d V_s=\int_0^t(\dot{\hbar}\circ U)(s)\d s, \ \ t\in[0,V_T]
\end{align*}
and
\beg{align*}
\int_0^{V_T}(\dot{\hbar}\circ U)^2(s)\d s=\int_0^T\dot{\hbar}^2(s)\d V_s<\infty.
\end{align*}

\emph{Proof of Proposition \ref{Pr1}.}
According to the linearity of $\mathcal{S}$, Lemma \ref{Le1} and the change of variables, we have for each $\hbar\in\mathfrak{H}_X$,
\beg{align}\label{Pf(Pr1)-1}
&\left(\mathcal{S}\left(\phi(X_t,\sL_{X_t})-\int_0^tb(s,X_s,\sL_{\vartheta(s,X_s,\sL_{X_s})})\d V_s\right)\right)(\hbar)\cr
&=(\mathcal{S}\phi(X_t,\sL_{X_t}))(\hbar)-\int_0^t(\mathcal{S}b(s,X_s,\sL_{\vartheta(s,X_s,\sL_{X_s})}))(\hbar)\d V_s\cr
&=(\tilde{\mathcal{S}}\phi(\tilde{W}_{V_t},\sL_{\tilde{W}_{V_t}}))(\hbar\circ U)-
\int_0^t(\tilde{\mathcal{S}}b(s,\tilde{W}_{V_s},\sL_{\vartheta(s,\tilde{W}_{V_s},\sL_{\tilde{W}_{V_s}})}))(\hbar\circ U)\d V_s\cr
&=(\tilde{\mathcal{S}}\phi(\tilde{W}_{V_t},\sL_{\tilde{W}_{V_t}}))(\hbar\circ U)-
\int_0^{V_t}(\tilde{\mathcal{S}}b(U_r,\tilde{W}_{r},\sL_{\vartheta(U_r,\tilde{W}_{r},\sL_{\tilde{W}_{r}})}))(\hbar\circ U)\d r\cr
&=\left(\tilde{\mathcal{S}}\left(\phi(\tilde{W}_{V_t},\sL_{\tilde{W}_{V_t}})-
\int_0^{V_t}b(U_r,\tilde{W}_{r},\sL_{\vartheta(U_r,\tilde{W}_{r},\sL_{\tilde{W}_{r}})})\d r\right)\right)(\hbar\circ U)\cr
&=\left(\tilde{\mathcal{S}}\int_0^{V_t}\si(U_r,\ti{W}_r,\sL_{\ti{W}_r})\d\ti{W}_r\right)(\hbar\circ U),
\end{align}
where the last equality is due to \eqref{Pr1-1}.\\
Since the classical It\^{o} integral coincides with the Wick-It\^{o} integral (see Remark \ref{Re(WI)} (ii))
and $\hbar\circ U$ belongs to $\mathfrak{H}_{\ti W}=\mathcal{R}(\mathfrak{C}_{\ti W})$, we get from \eqref{2Re(WI)}
\beg{align*}
\left(\tilde{\mathcal{S}}\int_0^{V_t}\si(U_r,\ti{W}_r,\sL_{\ti{W}_r})\d\ti{W}_r\right)(\hbar\circ U) =\int_0^{V_t}\left(\tilde{\mathcal{S}}\si(U_r,\ti{W}_r,\sL_{\ti{W}_r})\right)(\hbar\circ U)\d(\hbar\circ U)(r).
\end{align*}
Plugging this into \eqref{Pf(Pr1)-1} and using the change of variables and Lemma \ref{Le1} again, we obtain
\beg{align*}
&\left(\mathcal{S}\left(\phi(X_t,\sL_{X_t})-\int_0^tb(s,X_s,\P_{\vartheta(s,X_s,\sL_{X_s})})\d V_s\right)\right)(\hbar)\cr
&=\int_0^{V_t}\left(\tilde{\mathcal{S}}\si(U_r,\ti{W}_r,\sL_{\ti{W}_r})\right)(\hbar\circ U)\d(\hbar\circ U)(r)\cr
&=\int_0^t\left(\tilde{\mathcal{S}}\si(r,\ti{W}_{V_r},\sL_{\ti{W}_{V_r}})\right)(\hbar\circ U)\d\hbar(r)\cr
&=\int_0^t\left(\mathcal{S}\si(r,X_r,\sL_{X_r})\right)(\hbar)\d\hbar(r),
\end{align*}
which yields that $\int_0^t\si(s,X_s,\P_{X_s})\d^\diamond X_s$ is well-defined and moreover the relation \eqref{Pr1-2} holds due to Definition \ref{De-WI}.
The proof is now complete.
\qed

\subsection{The Lions derivative}

For later use, we state some basic facts about the Lions derivative.

For any $\th\in[1,\infty)$, $\sP_\th(\R^d)$ denotes the set of $\th$-integrable probability measures on $\R^d$,
and the $L^\th$-Wasserstein distance on $\sP_\th(\R^d)$ is defined as follows:
\begin{align*}
\mathbb{W}_\th(\mu,\nu):=\inf_{\pi\in\sC(\mu,\nu)}\left(\int_{\R^d\times\R^d}|x-y|^\th\pi(\d x, \d y)\right)^\ff 1 \th,\ \ \mu,\nu\in\sP_\th(\R^d),
\end{align*}
where $\sC(\mu,\nu)$ stands for the set of all probability measures on $\R^d\times\R^d$ with marginals $\mu$ and $\nu$.
It is well known that $(\sP_\th(\R^d),\mathbb{W}_\th)$ is a Polish space, usually referred to as the $\th$-Wasserstein space on $\R^d$.
We use $\<\cdot,\cdot\>$  for the Euclidean inner product, and $\|\cdot\|_{L^2_\mu}$ for the $L^2(\R^d\ra\R^d,\mu)$ norm.
Let $\sL_X$ be the distribution of  random variable $X$.

\beg{defn}
Let $f:\sP_2(\R^d)\ra\R$.
\begin{enumerate}
\item[(1)]  $f$ is called $L$-differentiable at $\mu\in\sP_2(\R^d)$, if the functional
\begin{align*}
L^2(\R^d\ra\R^d,\mu)\ni\phi\mapsto f(\mu\circ(\mathrm{Id}+\phi)^{-1}))
\end{align*}
is Fr\'{e}chet differentiable at $0\in L^2(\R^d\ra\R^d,\mu)$. That is, there exists a unique $\gamma\in L^2(\R^d\ra\R^d,\mu)$ such that
\begin{align*}
\lim_{\|\phi\|_{L^2_\mu}\ra0}\ff{f(\mu\circ(\mathrm{Id}+\phi)^{-1})-f(\mu)-\mu(\<\gamma,\phi\>)}{\|\phi\|_{L^2_\mu}}=0.
\end{align*}
In this case, $\gamma$ is called the $L$-derivative of $f$ at $\mu$ and denoted by $D^Lf(\mu)$.

\item[(2)] $f$ is called $L$-differentiable on $\sP_2(\R^d)$, if the $L$-derivative $D^Lf(\mu)$ exists for all $\mu\in\sP_2(\R^d)$.
Furthermore, if for every $\mu\in\sP_2(\R^d)$ there exists a $\mu$-version $D^Lf(\mu)(\cdot)$ such that $D^Lf(\mu)(x)$ is jointly continuous in $(\mu,x)\in\sP_2(\R^d)\times\R^d$, we denote $f\in C^{(1,0)}(\sP_2(\R^d))$.

\end{enumerate}
\end{defn}

In addition, according to \cite[Theorem 6.5]{Cardaliaguet13} and \cite[Proposition 3.1]{RW}, we get the following useful formula for the $L$-derivative.

\beg{lem}\label{FoLD}
Let $(\Omega,\sF,\P)$ be an atomless probability space and $\xi,\eta\in L^2(\Omega\ra\R^d,\P)$.
If $f\in C^{1,0}(\sP_2(\R^d))$, then
\begin{align*}
\lim_{\ve\da0}\ff {f(\sL_{\xi+\ve\eta})-f(\sL_\xi)} \ve=\E\<D^Lf(\sL_\xi)(\xi),\eta\>.
\end{align*}
\end{lem}

\section{Well-posedness of DDBSDE by Gaussian processes}

In this section, we consider the following DDBSDE driven by Gaussian process:
\begin{equation}\label{Bsde}
\left\{
\begin{array}{ll}
\d Y_t=-f(t,X_t,Y_t,Z_t,\sL_{(X_t,Y_t,Z_t)})\d V_t+Z_t\d^\diamond X_t,\\
Y_T=g(X_T,\sL_{X_T}),
\end{array} \right.
\end{equation}
where $f:[0,T]\times\R\times\R\times\R\times\sP_2(\R\times\R\times\R)\ra\R$ and $g:\R\times\sP_2(\R)\ra\R$ are measurable functions.
Now, let $\Upsilon$ be the space of pair of $(\sF_t^X)_{t\in[0,T]}$-adapted processes $(Y,Z)$ on $[0,T]$ satisfying
\begin{align*}
\E\left(\sup\limits_{0\leq t\leq T}|Y_t|^2+\int_0^T|Z_t|^2\d V_t\right)<\infty
\end{align*}
and
\beg{align*}
Y_t=u(t,X_t,\sL_{X_t}), \ \ \P{\textit-}a.s., \ t\in[0,T],\ \ Z_t=v(t,X_t,\sL_{X_t}), \ \d\P\otimes\d V{\textit-}a.e.
\end{align*}
with two deterministic functions $u,v:[0,T]\times\R\times\sP_2(\R)\ra\R$.

\beg{defn}\label{De-sol}
(1) A pair of stochastic processes $(Y,Z)=(Y_t,Z_t)_{0\leq t\leq T}$ is called a solution to \eqref{Bsde}, if $(Y,Z)\in\Upsilon$ satisfies
\beg{align*}
\E\left[g^2(X_T,\sL_{X_T})+\int_0^Tf^2(s,X_s,Y_s,Z_s,\sL_{(X_s,Y_s,Z_s)})\d V_s\right]<\infty,
\end{align*}
the Wick-It\^{o} integral $\int_0^tZ_s\d^\diamond X_s$ exists for any $t\in[0,T]$ and $\P$-a.s.
\beg{align*}
Y_t=g(X_T,\sL_{X_T})+\int_t^Tf(s,X_s,Y_s,Z_s,\sL_{(X_s,Y_s,Z_s)})\d V_s-\int_t^TZ_s\d^\diamond X_s,\ \ t\in[0,T].
\end{align*}

(2) \eqref{Bsde} is said to have uniqueness, if for any two solutions $(Y^i,Z^i)\in\Upsilon, i=1,2$,
\begin{align*}
\E\left(\sup\limits_{0\leq t\leq T}|Y^1_t-Y^2_t|^2+\int_0^T|Z^1_t-Z^2_t|^2\d V_t\right)=0.
\end{align*}

\end{defn}

In order to study \eqref{Bsde}, we shall introduce the following auxiliary BSDE
\begin{equation}\label{Aueq}
\left\{
\begin{array}{ll}
\d\tilde{Y}_t=-f(U_t,\tilde{W}_t,\tilde{Y}_t,\tilde{Z}_t,\sL_{(\tilde{W}_t,\tilde{Y}_t,\tilde{Z}_t)})\d t+\tilde{Z}_t\d\tilde{W}_t,\\
\tilde{Y}_{V_T}=g(\tilde{W}_{V_T},\sL_{\tilde{W}_{V_T}}),
\end{array} \right.
\end{equation}
where we recall that $\{\tilde{W}_t\}_{t\in[0,V_T]}$ is a standard Brownian motion stated in Section 2.2.
Let $\tilde{\Upsilon}$ be the space of pair of $(\sF_t^{\ti W})_{t\in[0,V_T]}$-adapted processes $(\ti Y,\ti Z)$ on $[0,V_T]$ satisfying
\begin{align*}
\ti\E\left(\sup\limits_{0\leq t\leq V_T}|\ti Y_t|^2+\int_0^{V_T}|\ti Z_t|^2\d t\right)<\infty
\end{align*}
and
\beg{align*}
\ti Y_t=\ti u(t,\ti W_t,\sL_{\ti W_t}), \ \ti\P{\textit-}a.s., \ t\in[0,V_T], \ \ \ti Z_t=\ti v(t,\ti W_t,\sL_{\ti W_t}), \ \d\ti\P\otimes\d t{\textit-}a.e.
\end{align*}
with two deterministic functions $\ti u,\ti v:[0,V_T]\times\R\times\sP_2(\R)\ra\R$.
Similar to Definition \ref{De-sol}, we can give the notion of a solution $(\ti Y,\ti Z)=(\ti Y_t,\ti Z_t)_{0\leq t\leq V_T}$ for \eqref{Aueq} in the space $\tilde{\Upsilon}$.

\beg{thm}\label{Th1}
Suppose that $\E [g^2(X_T,\sL_{X_T})]<\infty$ and  \eqref{Aueq} has a unique solution $(\tilde{Y},\tilde{Z})\in\tilde{\Upsilon}$.
Then  \eqref{Bsde} has a unique solution $(Y,Z)\in\Upsilon$.
More precisely, if
\beg{align*}
\tilde{Y}_t=\tilde{u}(t,\tilde{W}_t,\sL_{\tilde{W}_t}), \ \ \tilde{Z}_t=\tilde{v}(t,\tilde{W}_t,\sL_{\tilde{W}_t}),  \ \ t\in[0,V_T]
\end{align*}
is a unique solution of \eqref{Aueq},  then
\beg{align*}
Y_t=u(t,X_t,\sL_{X_t}), \ \ Z_t=v(t,X_t,\sL_{X_t}), \ \  t\in[0,T],
\end{align*}
is a unique solution of \eqref{Bsde}, where
\beg{align*}
u(t,x,\mu)=\tilde{u}(V_t,x,\mu),\ \ v(t,x,\mu)=\tilde{v}(V_t,x,\mu),\ \ (t,x,\mu)\in[0,T]\times\R\times\sP_2(\R).
\end{align*}
\end{thm}

\beg{proof}
The proof is divided into two steps.

\textit{Step 1. Existence.}
Since $(\tilde{Y},\tilde{Z})\in\tilde{\Upsilon}$ is a solution to \eqref{Aueq}, we let
$$\tilde{Y}_t=\tilde{u}(t,\tilde{W}_t,\sL_{\tilde{W}_t}),\ \ti\P{\textit-}a.s., \ t\in[0,V_T],  \ \
\tilde{Z}_t=\tilde{v}(t,\tilde{W}_t,\sL_{\tilde{W}_t}),\ \d\ti\P\otimes\d t{\textit-}a.e., $$
and then we have for any $t\in[0,T]$,
\beg{align*}
&\tilde{u}(V_t,\tilde{W}_{V_t},\sL_{\tilde{W}_{V_t}})-\tilde{u}(0,0,\delta_0)\cr
&=-\int_0^{V_t}f(U_s,\tilde{W}_s,\tilde{u}(s,\tilde{W}_{s},\sL_{\tilde{W}_{s}}),\tilde{v}(s,\tilde{W}_{s},\sL_{\tilde{W}_{s}}),
\sL_{(\tilde{W}_s,\tilde{u}(s,\tilde{W}_{s},\sL_{\tilde{W}_{s}}),\tilde{v}(s,\tilde{W}_{s},\sL_{\tilde{W}_{s}}))})\d s\cr
&\quad+\int_0^{V_t}\tilde{v}(s,\tilde{W}_{s},\sL_{\tilde{W}_{s}})\d\tilde{W}_s,
\end{align*}
where $\delta_0$ is the Dirac measure.\\
Because of the definition of $\tilde{\Upsilon}$,
it is readily checked that \eqref{Pr1-0} with $(\phi,b,\si)$ replaced by $(\tilde{u}(V_t,\cdot,\cdot),f,\tilde{v})$ holds.
So, by Proposition \ref{Pr1} we derive for any $t\in[0,T]$,
\beg{align*}
&\tilde{u}(V_t,X_t,\sL_{X_t})-\tilde{u}(0,0,\delta_0)\cr
&=-\int_0^tf(s,X_s,\tilde{u}(V_s,X_{s},\sL_{X_s}),\tilde{v}(V_s,X_s,\sL_{X_s}),
\sL_{(X_s,\tilde{u}(V_s,X_{s},\sL_{X_{s}}),\tilde{v}(V_s,X_{s},\sL_{X_{s}}))})\d V_s\cr
&\quad+\int_0^t\tilde{v}(V_s,X_s,\sL_{X_s})\d^\diamond X_s.
\end{align*}
Let $u(t,x,\mu)=\tilde{u}(V_t,x,\mu), v(t,x,\mu)=\tilde{v}(V_t,x,\mu),(t,x,\mu)\in[0,T]\times\R\times\sP_2(\R)$.
Then $(Y_t,Z_t):=(u(t,X_t,\sL_{X_t}),v(t,X_t,\sL_{X_t}))$ solves the following BSDE:
\beg{align}\label{Pf(Th1)-1}
Y_t=Y_T+\int_t^Tf(s,X_s,Y_s,Z_s,\sL_{(X_s,Y_s,Z_s)})\d V_s+\int_t^TZ_sd^\diamond X_s.
\end{align}
On the other hand, noting that
$g(\tilde{W}_{V_T},\sL_{\tilde{W}_{V_T}})=\tilde{u}(V_T,\tilde{W}_{V_T},\sL_{\tilde{W}_{V_T}}),\ \P{\textit-}a.s.$                                                                                                                  and the Gaussian law of $\tilde{W}_{V_T}$, we get
$$g(x,\sL_{\tilde{W}_{V_T}})=\tilde{u}(V_T,x,\sL_{\tilde{W}_{V_T}}),\ \ \d x{\textit-}a.e.,$$
which implies
\beg{align}\label{Pf(Th1)-2}
&g(X_T,\sL_{X_T})=g(X_T,\sL_{\tilde{W}_{V_T}})=\tilde{u}(V_T,X_T,\sL_{\tilde{W}_{V_T}})\cr
&=\tilde{u}(V_T,X_T,\sL_{X_T})=u(T,X_T,\sL_{X_T})=Y_T,\  \ \P{\textit-}a.s.
\end{align}
In addition, we easily obtain
\beg{align}\label{Pf(Th1)-3}
&\E\left(\sup\limits_{0\leq t\leq T}|Y_t|^2+\int_0^T|Z_t|^2\d V_t\right)\cr
&=\E\left(\sup\limits_{0\leq t\leq T}|u(t,X_t,\sL_{X_t})|^2+\int_0^T|v(t,X_t,\sL_{X_t})|^2\d V_t\right)\cr
&=\E\left(\sup\limits_{0\leq t\leq T}|\tilde{u}(V_t,X_t,\sL_{X_t})|^2+\int_0^T|\tilde{v}(V_t,X_t,\sL_{X_t})|^2\d V_t\right)\cr
&=\ti\E\left(\sup\limits_{0\leq t\leq T}|\tilde{u}(V_t,\tilde{W}_{V_t},\sL_{\tilde{W}_{V_t}})|^2+\int_0^T|\tilde{v}(V_t,\tilde{W}_{V_t},\sL_{\tilde{W}_{V_t}})|^2\d V_t\right)\cr
&=\ti\E\left(\sup\limits_{0\leq t\leq V_T}|\tilde{u}(t,\tilde{W}_{t},\sL_{\tilde{W}_{t}})|^2+\int_0^{V_T}|\tilde{v}(t,\tilde{W}_{t},\sL_{\tilde{W}_{t}})|^2\d t\right)\cr
&=\ti\E\left(\sup\limits_{0\leq t\leq V_T}|\tilde{Y}_t|^2+\int_0^{V_T}|\tilde{Z}_t|^2\d V_t\right)<\infty,
\end{align}
which, together with \eqref{Pf(Th1)-1}-\eqref{Pf(Th1)-2}, yields that $(Y,Z)\in\Upsilon$ is a solution to \eqref{Bsde}.

\textit{Step 2. Uniqueness.}
Let $(Y^i,Z^i)\in\Upsilon,i=1,2$ be two solutions of \eqref{Bsde}.
Then there exist $u^i,v^i,i=1,2$ such that $(Y_t^i,Z_t^i)=(u^i(t,X_t,\sL_{X_t}),v^i(t,X_t,\sL_{X_t})),i=1,2$.
Along the same arguments as in step 1, we can conclude that
$(\tilde{Y}_t^i,\tilde{Z}_t^i)=(u^i(U_t,\tilde{W}_t,\sL_{\tilde{W}_t}),v^i(U_t,\tilde{W}_t,\sL_{\tilde{W}_t}))\in\tilde{\Upsilon},i=1,2$ are two solutions of \eqref{Aueq}.
Similar to \eqref{Pf(Th1)-3} and by the uniqueness of \eqref{Aueq}, we arrive at
\beg{align*}
&\E\left(\sup\limits_{0\leq t\leq T}|Y_t^1-Y_t^2|^2+\int_0^T|Z_t^1-Z_t^2|^2\d V_t\right)\cr
&=\E\left(\sup\limits_{0\leq t\leq T}|u^1(t,X_t,\sL_{X_t})-u^2(t,X_t,\sL_{X_t})|^2+\int_0^T|v^1(t,X_t,\sL_{X_t})-v^2(t,X_t,\sL_{X_t})|^2\d V_t\right)\cr
&=\ti\E\left(\sup\limits_{0\leq t\leq V_T}|u^1(U_t,\tilde{W}_t,\sL_{\tilde{W}_t})-u^2(t,\tilde{W}_t,\sL_{\tilde{W}_t})|^2+\int_0^{V_T}|v^1(U_t,\tilde{W}_t,\sL_{\tilde{W}_t})-v^2(t,\tilde{W}_t,\sL_{\tilde{W}_t})|^2\d t\right)\cr
&=\ti\E\left(\sup\limits_{0\leq t\leq V_T}|\tilde{Y}_t^1-\tilde{Y}_t^2|^2+\int_0^{V_T}|\tilde{Z}_t^1-\tilde{Z}_t^2|^2\d t\right)=0,
\end{align*}
which finishes the proof.
\end{proof}

Next, we are devoted to applying our general Theorem \ref{Th1} to the case of Lipschitz continuous functions $(g,f)$.
More precisely, we assume the following conditions:
\begin{enumerate}

\item[\textsc{\textbf{(H1)}}] (i) $g,f$ are Lipschitz continuous,
i.e. there exist two constants $L_g, L_f>0$ such that
for all $t\in[0,T],x_i,y_i,z_i\in\R,\mu,\tilde{\mu}\in\sP_2(\R),\mu_i\in\sP_2(\R\times\R\times\R),i=1,2$,
\begin{align*}
|g(x_1,\mu)-g(x_2,\tilde{\mu})|\leq L_g(|x_1-x_2|+\mathbb{W}_2(\mu,\tilde{\mu}))
\end{align*}
and
\begin{align*}
|f(t,x_1,y_1,z_1,\mu_1)-f(t,x_2,y_2,z_2,\mu_2)|\leq L_f(|x_1-x_2|+|y_1-y_2|+|z_1-z_2|+\mathbb{W}_2(\mu_1,\mu_2)).
\end{align*}

\item[(ii)] $|g(0,\delta_0)|+\sup_{t\in [0,T]}|f(t,0,0,0,\delta_0)|<\infty$.
\end{enumerate}

Owing to \cite[Theorem A.1]{Li18} where the driven noises are a Brownian motion and an independent Poisson random measure,
the auxiliary equation \eqref{Aueq} has a unique solution $(\tilde{Y},\tilde{Z})\in\tilde{\Upsilon}$.
Hence, with the help of Theorem \ref{Th1}, \eqref{Bsde} admits a unique solution $(Y,Z)\in\Upsilon$.

\beg{rem}\label{Re1}

(i) (Stability estimate) Let $(Y^i,Z^i)\in\Upsilon$ be the unique solution of \eqref{Bsde} with the coefficients $(g^i,f^i),i=1,2$,
which satisfy \textsc{\textbf{(H1)}}. Then, there exists a constant $C>0$ such that
\begin{align*}
&\E\left(\sup_{s\in[0,T]}|Y^1_s-Y^2_s|^2+\int_0^{T}|Z^1_s-Z^2_s|^2\d V_s\right)\cr
\leq& C\E\left[\left|(g^1-g^2)(X_T,\sL_{X_{T}})\right|^2+
\left(\int_0^{T}|(f^1-f^2)(s,X_s,{Y}^1_s,{Z}^1_s,\sL_{(X_s,{Y}^1_s,{Z}^1_s)})|\d V_s\right)^2\right].
\end{align*}
Indeed, we denote by $(\ti Y^i,\ti Z^i)\in\ti\Upsilon$ the unique solution of \eqref{Aueq} with the coefficient $(g^i,f^i),i=1,2$.
By \cite[Theorem A.2]{Li18}, we derive that there is a constant $C>0$ such that
\begin{align*}
&\ti\E\left(\sup_{s\in[0,V_T]}|\ti Y^1_s-\ti Y^2_s|^2+\int_0^{V_T}|\ti Z^1_s-\ti Z^2_s|^2\d s\right)\cr
\leq& C\ti\E\left[\left|(g^1-g^2)(\tilde{W}_{V_T},\sL_{\tilde{W}_{V_T}})\right|^2+
\left(\int_0^{V_T}|(f^1-f^2)(U_s,\tilde{W}_s,\tilde{Y}^1_s,\tilde{Z}^1_s,\sL_{(\tilde{W}_s,\tilde{Y}^1_s,\tilde{Z}^1_s)})|\d s\right)^2\right].
\end{align*}
Then, the desired stability estimate follows from straightforward calculations via transfer principle Proposition \ref{Pr1}.

(ii) In the distribution-free case, namely the coefficients $g$ and $f$ in \eqref{Bsde} are independent of distributions,
Bender \cite[Theorems 4.2 and 4.4]{Bender14} proved the existence and uniqueness result, and our recent work \cite[Theorem 5.1]{FW21}
investigated the existence of densities and moreover derived their Gaussian estimates for the marginal laws of the solution.
So, our work in this note can be regarded as a continuation and generalization of \cite{Bender14,FW21}.
\end{rem}

\section{Comparison theorem and converse comparison theorem}

In this section, we first obtain a comparison theorem by imposing on the condition involved the Lions derivative of the generator.
In the second part, we concern with the converse problem of comparison theorem.
More precisely, we shall establish a representation theorem for the generator, and show how to combine this result to derive a converse comparison theorem.

\subsection{Comparison theorem}

In the Brownian motion case, i.e. $X=B^{1/2}$, the counter examples given in \cite[Examples 3.1 and 3.2]{BLP09} or \cite[Example 2.1]{LLZ18} show that
if the generator $f$ depends on the law of $Z$ or is not increasing with respect to the law of $Y$, comparison theorems usually don't hold for equation \eqref{Bsde}.
We consider now the following special version of (3.1):
\begin{align}\label{Bsde-com}
Y_t=g(X_T,\sL_{X_T})+\int_t^Tf(s,X_s,Y_s,Z_s,\sL_{(X_s,Y_s)})\d V_s-\int_t^TZ_s\d^\diamond X_s,\ \ t\in[0,T].
\end{align}
Correspondingly, the auxiliary equation of \eqref{Bsde-com} is of the form
\begin{align}\label{AuBsde-com}
\ti Y_t=&g(\tilde{W}_{V_T},\sL_{\ti{W}_{V_T}})+\int_t^{V_T}f(U_s,\ti{W}_s,\ti{Y}_s,\ti{Z}_s,\sL_{(\tilde{W}_s,\tilde{Y}_s)})\d s
-\int_t^{V_T}\tilde{Z}_s\d\tilde{W}_s,\ \ t\in[0,V_T].
\end{align}

\beg{thm}\label{Th(com)}
Assume that two generators $f^i,i=1,2$ satisfy \textsc{\textbf{(H1)}}, and that for any $(t,x,y,z)\in[0,T]\times\R\times\R\times\R$, $f^1(t,x,y,z,\cdot)$ belongs to $C^{(1,0)}(\sP_2(\R\times\R))$ with $0\leq (D^Lf^1)^{(2)}\leq K$ for some constant $K>0$, where $(D^Lf^1)^{(2)}$ denotes the second component of $D^Lf^1$.
Let $(Y^i,Z^i)\in\Upsilon, i=1,2$ be the solutions of \eqref{Bsde-com} with data $(g^i,f^i),i=1,2$, respectively.
Then, if
\begin{align}\label{1Th(com)}
g^1(x,\mu)\leq g^2(x,\mu),\ \ x\in\R,\mu\in\sP_2(\R)
\end{align}
and
\begin{align}\label{2Th(com)}
f^1(t,x,y,z,\nu)\leq f^2(t,x,y,z,\nu),\ \ t\in[0,T], x,y,z\in\R,\nu\in\sP_2(\R\times\R),
\end{align}
there holds that for every $t\in[0,T], Y_t^1\leq Y_t^2, \P$-a.s..
\end{thm}

\beg{proof}
Let $(\tilde{Y}^i,\tilde{Z}^i)\in\tilde{\Upsilon}$ be the unique solutions of \eqref{AuBsde-com} associated with $(g^i,f^i), i=1,2$,
and denote by $(\ti u^i,\ti v^i),i=1,2$ their representation functions, respectively.

We first give a comparison result for $(\tilde{Y}^i,\tilde{Z}^i), i=1,2$.
Our strategy hinges on the It\^{o}-Tanaka formula applied to $[(\ti Y_t^1-\ti Y_t^2)^+]^2$ (see, e,g., \cite[Proposition 2.35 and Remark 2.36]{PR14}), which gives for any $t\in[0,V_T]$,
\begin{align*}
&\ti\E[(\ti Y_t^1-\ti Y_t^2)^+]^2+\ti\E\int_t^{V_T}|\ti Z_s^1-\ti Z_s^2|^2\mathrm{I}_{\{\ti Y_s^1-\ti Y_s^2\geq0\}}\d s\nonumber\\
=&\ti\E[(g^1-g^2)^+(\tilde{W}_{V_T},\sL_{\ti{W}_{V_T}})]^2\nonumber\\
&+2\ti\E\int_t^{V_T}(\ti Y_s^1-\ti Y_s^2)^+\left(f^1(U_s,\Theta_s^1,\sL_{(\tilde{W}_s,\tilde{Y}^1_s)})-f^2(U_s,\Theta_s^2,\sL_{(\tilde{W}_s,\tilde{Y}^2_s)})\right)\d s\nonumber\\
=&2\ti\E\int_t^{V_T}(\ti Y_s^1-\ti Y_s^2)^+\left(f^1(U_s,\Theta_s^1,\sL_{(\tilde{W}_s,\tilde{Y}^1_s)})-f^2(U_s,\Theta_s^2,\sL_{(\tilde{W}_s,\tilde{Y}^2_s)})\right)\d s,
\end{align*}
Here, we have set $\Theta_s^i=(\tilde{W}_s,\ti{Y}^i_s,\ti{Z}^i_s), i=1,2$ and used \eqref{1Th(com)}.\\
Using the Lipschitz continuity of $f^1$ and the fact that $f^1(U_t,\cdot,\cdot,\cdot,\cdot)\leq f^2(U_t,\cdot,\cdot,\cdot,\cdot)$ for each $t\in[0,V_T]$ implied by \eqref{2Th(com)}, we get
\begin{align}\label{1PfTh(com)}
&\ti\E[(\ti Y_t^1-\ti Y_t^2)^+]^2+\ti\E\int_t^{V_T}|\ti Z_s^1-\ti Z_s^2|^2\mathrm{I}_{\{\ti Y_s^1-\ti Y_s^2\geq0\}}\d s\nonumber\\
=&2\ti\E\int_t^{V_T}(\ti Y_s^1-\ti Y_s^2)^+
\left(f^1(U_s,\Theta_s^1,\sL_{(\tilde{W}_s,\tilde{Y}^1_s)})-f^1(U_s,\Theta_s^2,\sL_{(\tilde{W}_s,\tilde{Y}^1_s)})\right)\d s\cr
&+2\ti\E\int_t^{V_T}(\ti Y_s^1-\ti Y_s^2)^+
\left(f^1(U_s,\Theta_s^2,\sL_{(\tilde{W}_s,\tilde{Y}^1_s)})-f^1(U_s,\Theta_s^2,\sL_{(\tilde{W}_s,\tilde{Y}^2_s)})\right)\d s\cr
&+2\ti\E\int_t^{V_T}(\ti Y_s^1-\ti Y_s^2)^+
\left(f^1(U_s,\Theta_s^2,\sL_{(\tilde{W}_s,\tilde{Y}^2_s)})-f^2(U_s,\Theta_s^2,\sL_{(\tilde{W}_s,\tilde{Y}^2_s)})\right)\d s\cr
\leq&2L_{f^1}\ti\E\int_t^{V_T}(\ti Y_s^1-\ti Y_s^2)^+\left(|\ti Y_s^1-\ti Y_s^2|+|\ti Z_s^1-\ti Z_s^2|\right)\d s
\cr
&+2\ti\E\int_t^{V_T}(\ti Y_s^1-\ti Y_s^2)^+
\left(f^1(U_s,\Theta_s^2,\sL_{(\tilde{W}_s,\tilde{Y}^1_s)})-f^1(U_s,\Theta_s^2,\sL_{(\tilde{W}_s,\tilde{Y}^2_s)})\right)\d s.
\end{align}
Observe that by Lemma \ref{FoLD}, we obtain
\begin{align*}
&f^1(U_s,\Theta_s^2,\sL_{(\tilde{W}_s,\tilde{Y}^1_s)})-f^1(U_s,\Theta_s^2,\sL_{(\tilde{W}_s,\tilde{Y}^2_s)})\cr
=&\int_0^1\ff \d {\d r}f^1(U_s,\Theta_r^2,\sL_{\chi_s(r)})\d r\cr
=&\int_0^1\ti\E\langle D^Lf^1(U_s,x,\cdot)(\sL_{\chi_s(r)})(\chi_s(r)),(0,\tilde{Y}^1_s-\tilde{Y}^2_s)\rangle|_{x=\Theta_r^2}\d r\cr
=&\int_0^1\ti\E\left[(D^Lf^1)^{(2)}(U_s,x,\cdot)(\sL_{\chi_s(r)})(\chi_s(r))\cdot(\tilde{Y}^1_s-\tilde{Y}^2_s)\right]\big|_{x=\Theta_r^2}\d r\cr
\leq& K\ti\E(\tilde{Y}^1_s-\tilde{Y}^2_s)^+,
\end{align*}
where for any $r\in[0,1],\chi_s(r)=(\tilde{W}_s,\tilde{Y}^2_s)+r(0,\tilde{Y}^1_s-\tilde{Y}^2_s)$.\\
Therefore, plugging this into \eqref{1PfTh(com)} and applying the H\"{o}lder and the Young inequalities, we have
\begin{align*}
&\ti\E[(\ti Y_t^1-\ti Y_t^2)^+]^2+\ti\E\int_t^{V_T}|\ti Z_s^1-\ti Z_s^2|^2\mathrm{I}_{\{\ti Y_s^1-\ti Y_s^2\geq0\}}\d s\nonumber\\
\leq&2(L_{f^1}\vee K)\ti\E\int_t^{V_T}(\ti Y_s^1-\ti Y_s^2)^+\left(|\ti Y_s^1-\ti Y_s^2|+|\ti Z_s^1-\ti Z_s^2|+\ti\E(\tilde{Y}^1_s-\tilde{Y}^2_s)^+\right)\d s\cr
\leq&C\int_t^{V_T}\ti\E[(\ti Y_s^1-\ti Y_s^2)^+]^2\d s+\ff 1 2\ti\E\int_t^{V_T}|\ti Z_s^1-\ti Z_s^2|^2\mathrm{I}_{\{\ti Y_s^1-\ti Y_s^2\geq0\}}\d s,
\end{align*}
where and in what follows C denotes a generic constant.\\
Then, the Gronwall inequality implies $\ti\E[(\ti Y_t^1-\ti Y_t^2)^+]^2=0$ for every $t\in[0,V_T]$.
Consequently, it holds that for any $t\in[0,V_T]$,
\begin{align*}
\tilde{u}^1(t,\tilde{W}_t,\sL_{\tilde{W}_t})=\ti Y_t^1\leq \ti Y_t^2=\tilde{u}^2(t,\tilde{W}_t,\sL_{\tilde{W}_t}),\ \  \ti\P\textit{-}a.s..
\end{align*}
Since $\tilde{W}_t$ has the Gaussian law, we derive that for any $t\in[0,V_T]$,
\begin{align*}
\tilde{u}^1(t,x,\sL_{\tilde{W}_t})\leq\tilde{u}^2(t,x,\sL_{\tilde{W}_t}),\ \  \d x\textit{-}a.e.,
\end{align*}
which, together with the fact that $\sL_{\tilde{W}_{V_t}}=\sL_{X_t}$, yields that for each $t\in[0,T]$,
\begin{align*}
\tilde{u}^1(V_t,x,\sL_{X_t})\leq\tilde{u}^2(V_t,x,\sL_{X_t}),\ \  \d x\textit{-}a.e.,
\end{align*}
According to Theorem \ref{Th1}, we conclude that $Y_t^1\leq Y_t^2, \P$-a.s. for every $t\in[0,T]$.
This completes the proof.
\end{proof}

\beg{rem}\label{Re-comp1}
(i) In light of the proof above, one can see that if the conditions for $f^1$ are replaced by that for $f^2$,
Theorem \ref{Th(com)} also holds.

(ii) Compared with the relevant result on DDBSDEs driven by the standard Brownian motion $B^{1/2}$  shown in \cite[Theorem 2.2 and Remark 2.2]{LLZ18},
it is easy to see that our above result Theorem \ref{Th(com)} applies to more general BSDEs since we replace $B^{1/2}$ with a general
Gaussian process $X$ as driving process.
In addition, in contrast to the distribution-free case (see, e.g.,  \cite[Theorem 4.6 (iii)]{Bender14}, \cite[Theorem 2.2]{EPQ9705} or \cite[Theorem 12.3]{HCS12}), we need to overcome difficulty induced by the appearance of $\sL_{(X_s,Y_s)}$ in the generator via a formula for the $L$-derivative (Lemma \ref{FoLD}).
\end{rem}

\subsection{Converse comparison theorem}

In this part, we investigate a kind of converse comparison problem:
if for each $t\in[0,T],\ep\in(0,T-t)$ and $\xi\in L^2_X$, $Y_t^1(t+\ep,\xi)\leq Y_t^2(t+\ep,\xi)$ (see the definition at the beginning of Theorem \ref{Th(Conve)}), do we have $f^1(t,y,z,\nu)\leq f^2(t,y,z,\nu)$ for every $(t,y,z,\nu)\in[0,T]\times\R\times\R\times\sP(\R\times\R\times\R)$?
That is, if we can compare the solutions of two DDBSDEs with the
same terminal condition, for all terminal conditions, can we compare the generators?
In order to study this question, we first give a representation theorem for the generator $f$ initiated in \cite{BCHMP00}.

Now, given $(t,y,z)\in[0,T)\times\R\times\R$ and let $0<\epsilon<T-t$, and denote by $(Y^\epsilon,Z^\epsilon)$ the unique solution of the following DDBSDE on $[t,t+\epsilon]$:
\begin{align}\label{1-Repr}
Y^\epsilon_r=y+z(X_{t+\ep}-X_t)+\int_r^{t+\ep}f(s,Y^\ep_s,Z^\ep_s,\sL_{(X_s,Y_s^\ep,Z_s^\ep)})\d V_s-\int_r^{t+\epsilon}Z^\epsilon_s\d^\diamond X_s.
\end{align}

Our representation theorem is formulated as follows

\beg{thm}\label{Th(Repre)}
Assume that \textsc{\textbf{(H1)}} holds and $V_{t+\ep}-V_t=O(\ep)$ as $\ep\ra0$ for any $t\in[0,T)$.
Then for any $(t,y,z)\in[0,T)\times\R\times\R$, the following two statements are equivalent:
\beg{align}\label{1Th(Re)}
(i) \lim\limits_{\ep\ra0^+}&\frac{Y^\ep_t-y}{\ep}=f(t,y,z,\sL_{(X_t,y,z)});\\
(ii) \lim\limits_{\ep\ra0^+}&\ff 1 \ep\int_t^{t+\ep}f(r,y,z,\sL_{(X_t,y,z)})\d r=f(t,y,z,\sL_{(X_t,y,z)}).\label{2Th(Re)}\ \ \  \ \ \  \ \ \  \ \ \  \ \ \ \ \ \  \ \ \  \  \  \  \  \ \  \  \ \  \  \ \
\end{align}
\end{thm}

\beg{proof}
We split the proof into two steps.
First, we prove this theorem for the case of the auxiliary equation of \eqref{1-Repr};
then we extend this result to \eqref{1-Repr} by applying the transfer principle.

\textit{Step 1. The auxiliary equation case.}
Observe first that the corresponding auxiliary equation of \eqref{1-Repr} is given by
\begin{align}\label{1PfTh(Re)}
\ti{Y}^\ep_{V_r}=&y+z(\ti{W}_{V_{t+\ep}}-\ti{W}_{V_t})+\int_{V_r}^{V_{t+\ep}}f(U_s,\ti{Y}^\ep_s,\ti{Z}^\ep_s,\sL_{(\ti{W}_s,\ti{Y}^\ep_s,\ti{Z}_s^\ep)})\d s\cr
&-\int_{V_r}^{V_{t+\ep}}\ti{Z}^\ep_s\d\ti{W}_s,\ \ r\in[t,t+\ep].
\end{align}
We claim that for any $(t,y,z)\in[0,T)\times\R\times\R$, the following two statements are equivalent:
\begin{align*}
(I) \lim\limits_{\ep\ra0^+}&\frac{\ti Y^\epsilon_{V_t}-y}{\epsilon}=f(t,y,z,\sL_{(\ti W_{V_t}.y,z)});\\
(II) \lim\limits_{\ep\ra0^+}&\ff 1 \ep\int_{V_t}^{V_{t+\ep}}f(U_r,y,z,\sL_{(\ti W_{V_t}.y,z)})\d r=f(t,y,z,\sL_{(\ti W_{V_t}.y,z)}).
\ \ \  \ \ \  \ \ \  \ \ \  \ \ \ \ \ \  \   \ \  \  \ \  \  \ \  \  \ \  \  \ \  \  \ \  \ \  \ \  \
\end{align*}
Indeed, for $s\in[V_t,V_{t+\ep}]$, we put $\Gamma^\ep_s:=\ti{Y}^\ep_s-(y+z(\ti{W}_s-\ti{W}_{V_t}))$.
Then, applying the It\^{o} formula to $\Gamma^\ep_s$ on the interval $[V_t,V_{t+\ep}]$  yields
\begin{align}\label{2PfTh(Re)}
\Gamma^\ep_s=&\int_s^{V_{t+\ep}}f(U_r,\Gamma^\ep_r+y+z(\ti{W}_r-\ti{W}_{V_t}),\ti{Z}^\ep_r,\sL_{(\ti{W}_r,\Gamma^\ep_r+y+z(\ti{W}_r-\ti{W}_{V_t}),\ti{Z}_r^\ep)})\d r\cr
&-\int_s^{V_{t+\ep}}(\ti{Z}^\ep_r-z)\d\ti{W}_r\cr
=:&\int_s^{V_{t+\ep}}f^\ep(r)\d r-\int_s^{V_{t+\ep}}(\ti{Z}^\ep_r-z)\d\ti{W}_r.
\end{align}
Using the facts that $\ti{Y}^\ep_{V_t}-y=\Gamma^\ep_{V_t}$ and $\Gamma^\ep_{V_t}$ is deterministic thanks to \cite[Proposition 4.2]{EPQ9705} and taking the conditional expectation with respect to $\ti\sF_{V_t}$ in \eqref{2PfTh(Re)}, we arrive at
\begin{align}\label{6PfTh(Re)}
&\frac{\ti Y^\epsilon_{V_t}-y}{\epsilon}-f(t,y,z,\sL_{(\ti W_{V_t}.y,z)})\cr
=&\ff{\Gamma^\ep_{V_t}}\ep-f(t,y,z,\sL_{(\ti W_{V_t}.y,z)})\cr
=&\ff 1 \ep\ti\E\left(\int_{V_t}^{V_{t+\ep}}f^\ep(r)\d r|\ti\sF_{V_t}\right)-f(t,y,z,\sL_{(\ti W_{V_t}.y,z)})\cr
=&\ff 1 \ep\ti\E\left[\int_{V_t}^{V_{t+\ep}}\left(f^\ep(r)-f(U_r,y+z(\ti{W}_r-\ti{W}_{V_t}),z,\sL_{(\ti{W}_r,y+z(\ti{W}_r-\ti{W}_{V_t}),z)})\right)\d r\big|\ti\sF_{V_t}\right]\cr
&+\ff 1 \ep\ti\E\bigg[\int_{V_t}^{V_{t+\ep}}\Big(f(U_r,y+z(\ti{W}_r-\ti{W}_{V_t}),z,\sL_{(\ti{W}_r,y+z(\ti{W}_r-\ti{W}_{V_t}),z)})\cr
&\qquad\qquad\qquad\quad -f(U_r,y,z,\sL_{(\ti W_{V_t}.y,z)})\Big)\d r\big|\ti\sF_{V_t}\bigg]\cr
&+\ff 1 \ep\int_{V_t}^{V_{t+\ep}}f(U_r,y,z,\sL_{(\ti W_{V_t}.y,z)})\d r-f(t,y,z,\sL_{(\ti W_{V_t}.y,z)})\cr
=:&I_1+I_2+I_3.
\end{align}
By the H\"{o}lder inequality and \textsc{\textbf{(H1)}}, we have
\begin{align}\label{3PfTh(Re)}
\ti\E|I_1|^2\leq&3L_f^2\ff {V_{t+\ep}-V_t}{\ep^2}\int_{V_t}^{V_{t+\ep}}\ti\E\left[|\Gamma^\ep_r|^2+|\ti{Z}^\ep_r-z|^2+\ti\E\left(|\Gamma^\ep_r|^2+|\ti{Z}^\ep_r-z|^2\right)\right]\d r\cr
=&6L_f^2\ff {V_{t+\ep}-V_t}{\ep^2}\int_{V_t}^{V_{t+\ep}}\ti\E\left(|\Gamma^\ep_r|^2+|\ti{Z}^\ep_r-z|^2\right)\d r
\end{align}
and
\begin{align}\label{4PfTh(Re)}
\ti\E|I_2|^2\leq&2L_f^2\ff {V_{t+\ep}-V_t}{\ep^2}\int_{V_t}^{V_{t+\ep}}\E\left[z^2|\ti{W}_r-\ti{W}_{V_t}|^2+(1+z^2)\ti\E|\ti{W}_r-\ti{W}_{V_t}|^2\right]\d r\cr
=&2L_f^2(1+2z^2)\ff {V_{t+\ep}-V_t}{\ep^2}\int_{V_t}^{V_{t+\ep}}(r-V_t)\d r\cr
=&L_f^2(1+2z^2)\ff {(V_{t+\ep}-V_t)^3}{\ep^2}.
\end{align}
Similar to \cite[Proposition 2.2]{BCHMP00}, using the It\^{o} formula applied to $\e^{\be s}|\Gamma^\ep_s|^2$ with some constant $\beta>0$ depending only on $L_f$, we obtain the following a priori estimate for \eqref{2PfTh(Re)} (its solution is regarded as $(\Gamma^\ep_\cdot,\ti{Z}^\ep_\cdot-z)$)
\begin{align*}
&\ti\E\left(\sup_{V_t\leq s\leq V_{t+\ep}}|\Gamma^\ep_s|^2+\int_{V_t}^{V_{t+\ep}}|\ti{Z}^\ep_r-z|^2\d r\big|\ti\sF_{V_t}\right)\cr
\leq&C\ti\E\left[\left(\int_{V_t}^{V_{t+\ep}}|
f(U_r,y+z(\ti{W}_r-\ti{W}_{V_t}),z,\sL_{(\ti{W}_r,y+z(\ti{W}_r-\ti{W}_{V_t}),z)})|\d r\right)^2\big|\ti\sF_{V_t}\right].
\end{align*}
Then, it follows from the H\"{o}lder inequality and \textsc{\textbf{(H1)}} that
\begin{align*}
&\ti\E\left(\sup_{V_t\leq s\leq V_{t+\ep}}|\Gamma^\ep_s|^2+\int_{V_t}^{V_{t+\ep}}|\ti{Z}^\ep_r-z|^2\d r\right)\cr
\leq&C(V_{t+\ep}-V_t)
\int_{V_t}^{V_{t+\ep}}\ti\E\left[|f(U_r,y,z,\sL_{(\ti{W}_{V_t},y,z)})|^2+z^2|\ti{W}_r-\ti{W}_{V_t}|^2+(1+z^2)\ti\E|\ti{W}_r-\ti{W}_{V_t}|^2\right]\d r\cr
=&C(V_{t+\ep}-V_t)
\left[\int_{V_t}^{V_{t+\ep}}|f(U_r,y,z,\sL_{(\ti{W}_{V_t},y,z)})|^2\d r+(1+2z^2)(V_{t+\ep}-V_t)^2\right].
\end{align*}
Substituting this into \eqref{3PfTh(Re)} yields
\begin{align}\label{5PfTh(Re)}
\ti\E|I_1|^2\leq& C\left(\ff{V_{t+\ep}-V_t}{\ep}\right)^2(1+V_{t+\ep}-V_t)\cr
&\times\left[\int_{V_t}^{V_{t+\ep}}|f(U_r,y,z,\sL_{(\ti{W}_{V_t},y,z)})|^2\d r+(1+2z^2)(V_{t+\ep}-V_t)^2\right].
\end{align}
Note that by the absolute continuity of integral, we have
\begin{align*}
\lim\limits_{\ep\ra0^+}\int_{V_t}^{V_{t+\ep}}|f(U_r,y,z,\sL_{(\ti{W}_{V_t},y,z)})|^2\d r=0.
\end{align*}
Therefore, owing to the condition $V_{t+\ep}-V_t=O(\ep)$ as $\ep\ra0$ and \eqref{4PfTh(Re)}-\eqref{5PfTh(Re)}, we conclude that
\begin{align*}
\lim\limits_{\ep\ra0^+}\left(\ti\E|I_1|^2+\ti\E|I_2|^2\right)=0,
\end{align*}
which, along with \eqref{6PfTh(Re)} and the fact that $\ti Y^\epsilon_{V_t}$ is deterministic due to \cite[Proposition 4.2]{EPQ9705} again, implies the desired assertion.

\textit{Step 2. The equation with Gaussian noise case.}
 With a slight modification of the proof of Theorem \ref{Th1}, we know that for every $\ep>0$, there exists a representation function $\ti u^\ep$ such that the solutions of \eqref{1-Repr} and \eqref{1PfTh(Re)} can be written as the following form: for any $r\in[t,t+\ep]$,
\begin{align*}
Y^\ep_r=\ti u^\ep(V_r,X_r-X_t,\sL_{X_r-X_t}) \ \ \ \mathrm{and} \ \ \ \ti Y^\ep_{V_r}=\ti u^\ep(V_r,\ti W_{V_r}-\ti W_{V_t},\sL_{\ti W_{V_r}-\ti W_{V_t}}).
\end{align*}
Then, we get
\begin{align}\label{7PfTh(Re)}
\ff {Y^\ep_t-y}\ep=\ff{\ti u^\ep(V_t,0,\de_0)-y}\ep=\ff{\ti Y^\ep_{V_t}-y}{\ep}.
\end{align}
On the other hand, by a change of variables it is easy to see that
\begin{align}\label{8PfTh(Re)}
\int_{t}^{t+\ep}f(r,y,z,\sL_{(\ti W_{V_t},y,z)})\d r=\int_{V_t}^{V_{t+\ep}}f(U_r,y,z,\sL_{(\ti W_{V_t}.y,z)})\d r.
\end{align}
Finally, taking into account of the claim in step 1 and using \eqref{7PfTh(Re)}-\eqref{8PfTh(Re)} and the fact that $\sL_{(\ti W_{V_t},y,z)}=\sL_{(X_t,y,z)}$,
we deduce that \eqref{1Th(Re)} is equivalent to \eqref{2Th(Re)}.
Our proof is now finished.
\end{proof}

\beg{rem}\label{Re(Rep)}
Although the computations become much more involved, it is not hard to extend the result in step 1 to \eqref{1PfTh(Re)} with $\ti W$ being replaced by a diffusion process, which is a generalisation of \cite[Proposition 2.3]{BCHMP00} and \cite[Theorem 3.3]{Jiang05} that handled the distribution-free BSDEs
driven by Brownian motion.
\end{rem}

With the help of Theorem \ref{Th(Repre)}, we can establish a converse comparison theorem.
To this end, we denote by $(Y_t^i(T,g(X_T,\sL_{X_T})),Z_t^i(T,g(X_T,\sL_{X_T})))_{t\in[0,T]}$ the solution of
\beg{align}\label{3Th(Conve)}
Y_t^i=&g(X_T,\sL_{X_T})+\int_t^Tf^i(s,Y_s^i,Z_s^i,\sL_{(X_s,Y^i_s,Z^i_s)})\d V_s\cr
&-\int_t^TZ_s^i\d^\diamond X_s,\ \ t\in[0,T],\ \ i=1,2.
\end{align}

Our converse comparison theorem reads as follows

\beg{thm}\label{Th(Conve)}
Let \textsc{\textbf{(H1)}} hold for $f^i,i=1,2$ and $V_{t+\ep}-V_t=O(\ep)$ as $\ep\ra0$ for any $t\in[0,T)$.
Assume moreover that for each $(t,y,z)\in[0,T)\times\R\times\R$ and $\ep\in(0,T-t]$,
\beg{align}\label{1Th(Conve)}
Y_t^1(t+\ep,y+z(X_{t+\ep}-X_t))\leq Y_t^2(t+\ep,y+z(X_{t+\ep}-X_t)).
\end{align}
Then for  each $(t,y,z)\in[0,T)\times\R\times\R$, we have
\beg{align}\label{2Th(Conve)}
f^1(t,y,z,\sL_{(X_t,y,z)})\leq f^2(t,y,z,\sL_{(X_t,y,z)}).
\end{align}
\end{thm}

\beg{proof}
Since \textsc{\textbf{(H1)}} holds for $f^i,i=1,2$, one can see that  for any $(t,y,z)\in[0,T)\times\R\times\R$,
\beg{align*}
\lim\limits_{\ep\ra0^+}&\ff 1 \ep\int_t^{t+\ep}f^i(r,y,z,\sL_{(X_t,y,z)})\d r=f^i(t,y,z,\sL_{(X_t,y,z)}),\ \ i=1,2.
\end{align*}
Then, in light of Theorem \ref{Th(Repre)}, we have
\beg{align}\label{1PfTh(Conve)}
\lim\limits_{\ep\ra0^+}&\frac{Y_t^i(t+\ep,y+z(X_{t+\ep}-X_t))-y}{\ep}=f^i(t,y,z,\sL_{(X_t,y,z)}),\ \ i=1,2.
\end{align}
By the hypothesis \eqref{1Th(Conve)} we derive that for any $(t,y,z)\in[0,T)\times\R\times\R$ and $\ep\in(0,T-t]$,
\beg{align*}
Y_t^1(t+\ep,y+z(X_{t+\ep}-X_t))-y\leq Y_t^2(t+\ep,y+z(X_{t+\ep}-X_t))-y,
\end{align*}
which, along with \eqref{1PfTh(Conve)}, yields the desired relation.
This then completes the proof.
\end{proof}

We conclude this section with a remark.

\beg{rem}\label{Re(Conve)}
By Theorems \ref{Th(Repre)} and \ref{Th(Conve)}, it is a surprise to us for finding that the relations \eqref{1Th(Re)}, \eqref{2Th(Re)} and \eqref{2Th(Conve)} all depend on $\sL_{(X_t,y,z)}$, which has a big difference from distribution-free case (see, e.g.,  \cite[Theorem 4.1]{BCHMP00}, \cite[Theorem]{Chen98a} or
\cite[Theroem 5.1]{Jiang05}).
This means that there is no arbitrariness for the measure of the generator when dealing with representation theorem or converse comparison theorem,
which is due to the appearance of the distribution dependent terms $\sL_{(X_s,Y_s^\ep,Z_s^\ep)}$ in \eqref{1-Repr} and $\sL_{(X_s,Y^i_s,Z^i_s)}$ in \eqref{3Th(Conve)}, respectively.
\end{rem}

\section{Functional inequalities}

In this section, we aim to establish functional inequalities for \eqref{Bsde}, including mainly transportation inequalities and Logarithmic-Sobolev inequalities.
Our arguments consist of utilising stability of the Wasserstein distance and the relative entropy of measures under the homeomorphism, together with the transfer principle.

\subsection{Transportation inequalities}

Let $(E,d)$ be a metric space and let $\sP(E)$ denote the set of all probability measures on $E$.
For $p\in[1,\infty)$, we say that a probability measure $\mu\in\sP(E)$ satisfies $p$-transportation inequality on $(E,d)$ (noted $\mu\in T_p(C)$) if there is a constant $C\geq0$ such that
\begin{align*}
\W_p(\mu,\nu)\leq C\sqrt{H(\nu|\mu)}, \ \ \nu\in\sP(E),
\end{align*}
where $H(\nu|\mu)$ is the relative entropy (or Kullback-Leibler divergence) of $\nu$ with respect to $\mu$ defined as
\begin{equation*}
H(\nu|\mu)=\left\{
\begin{array}{ll}
\int\log\ff {\d\nu}{\d\mu}\d\nu,\ \ \mathrm{if}\ \nu\ll\mu,\\
+\infty,\ \ \ \ \ \ \ \ \ \ \mathrm{otherwise}.
\end{array} \right.
\end{equation*}
The transportation inequality has found numerous applications, for instance, to quantitative finance, the concentration of measure phenomenon and various problems of probability in higher dimensions.
We refer the reader e.g. to \cite{BT20,DGW04,Lacker18,Rid17,Saussereau12,SYZ22} and references therein.

Before moving to \eqref{Bsde}, we first show the transportation inequalities for the auxiliary equation \eqref{Aueq}, which is a distribution dependent version of \cite[Theorem 1.3 and Lemma 4.1]{BT20}.

\beg{prp}\label{Prp(Tr)}
Assume that \textsc{\textbf{(H1)}} holds. Then we have
\begin{itemize}
\item[(i)] The law of $(\ti Y_t)_{t\in[0,V_T]}$ satisfies $T_2(C_{Tr,\ti Y})$ on $\ti\Om$ with $C_{Tr,\ti Y}=2(L_g+L_fV_T)^2\e^{2L_fV_T}$.
\item[(ii)] For any $p\geq1$, the law of $(\ti Z_t)_{t\in[0,V_T]}$, denoted by $\sL_{\ti Z}$, satisfies
\begin{align*}
\W_p(\sL_{\ti Z},\mu)\leq C_{Tr,\ti Z}\left(H(\mu|\sL_{\ti Z})\right)^{\ff 1 {2p}},\ \ \mu\in\sP(\ti\Om),
\end{align*}
where
\begin{align*}
C_{Tr,\ti Z}=2\inf_{\al>0}\left\{\ff 1{2\al}\left[1+\al\e^{2pL_fV_T}(L_g+L_fV_T)^{2p}\right]\right\}^{\ff 1{2p}}.
\end{align*}
\end{itemize}
\end{prp}

For the sake of conciseness, we defer the proof to the Appendix.
With Proposition \ref{Prp(Tr)} in hand, along with the transfer principle, we now state and prove the transportation inequalities for \eqref{Bsde} as follows.

\beg{thm}\label{Th(TrIn)}
Assume that \textsc{\textbf{(H1)}} holds. Then for every $t\in[0,T]$, we have
\begin{itemize}
\item[(i)] The law of $Y_t$ satisfies $T_2(C_{Tr,Y_t})$ on $\R$ with $C_{Tr,Y_t}=2(L_g+L_f(V_T-V_t))^2\e^{2L_f(V_T-V_t)}$.
\item[(ii)] For any $p\geq1$, the law of $Z_t$ satisfies
\begin{align*}
\W_p(\sL_{Z_t},\mu)\leq C_{Tr,Z_t}\left(H(\mu|\sL_{Z_t})\right)^{\ff 1 {2p}},\ \ \mu\in\sP(\R),
\end{align*}
where
\begin{align*}
C_{Tr,Z_t}=2\inf_{\al>0}\left\{\ff 1{2\al}\left[1+\al\e^{2pL_f(V_T-V_t)}(L_g+L_f(V_T-V_t))^{2p}\right]\right\}^{\ff 1{2p}}.
\end{align*}
\end{itemize}
\end{thm}

\beg{proof}
By Proposition \ref{Prp(Tr)} and its proof (see \eqref{add4PfPrp(Tr)} and \eqref{6PfPrp(Tr)} in the Appendix), it is easy to see that for each $t\in[0,V_T]$,
$\sL_{\ti Y_t}$ satisfies $T_2(C_{Tr,\ti Y_t})$ on $\R$ with $C_{Tr,\ti Y_t}=2(L_g+L_f(V_T-t))^2\e^{2L_f(V_T-t)}$, and $\sL_{\ti Z_t}$ satisfies
\begin{align*}
\W_p(\sL_{\ti Z_t},\mu)\leq C_{Tr,\ti Z_t}\left(H(\mu|\sL_{\ti Z_t})\right)^{\ff 1 {2p}},\ \ \mu\in\sP(\R)
\end{align*}
with any $p\geq1$ and
\begin{align*}
C_{Tr,\ti Z_t}=2\inf_{\al>0}\left\{\ff 1{2\al}\left[1+\al\e^{2pL_f(V_T-t)}(L_g+L_f(V_T-t))^{2p}\right]\right\}^{\ff 1{2p}}.
\end{align*}
Noting that for any $t\in[0,T]$,  the laws of $Y_t$ and $Z_t$ are respectively the same as those of $\ti Y_{V_t}$ and $\ti Z_{V_t}$ due to Theorem \ref{Th1},
we obtain the desired assertions (i) and (ii).
\end{proof}

\beg{rem}\label{Re(TrIn)}
If $f=0$ and $g(x,\mu)=x$, then one can check that \eqref{Bsde} and \eqref{Aueq} have unique solutions $(Y,Z)=(X,1)$ and $(\ti Y,\ti Z)=(\ti W,1)$, respectively.
By Theorem \ref{Th(TrIn)} and Proposition \ref{Prp(Tr)}, we have $C_{Tr,Y_t}=2$ and $C_{Tr,\ti Y}=2$ which are known to be optimal for Gaussian processes and Brownian motion.
This means that the constants $C_{Tr,Y_t}$ and $C_{Tr,\ti Y}$ above are sharp.
Let us also mention that it is not clear here, if the laws of the paths of $Y$ and  $Z$ satisfy the transportation inequality or not.
Indeed, the transfer principle may fail in this situation because of the difference between spaces $\sP(C([0,T]))$ and $\sP(C([0,V_T]))$.
\end{rem}

\subsection{Logarithmic-Sobolev inequality}

For introducing our result for \eqref{Bsde}, let us first give the Logarithmic-Sobolev inequality for the auxiliary equation \eqref{Aueq}.

\beg{prp}\label{Prp(LS)}
Assume that \textsc{\textbf{(H1)}} holds. Then for every $t\in[0,V_T]$,
\begin{align}\label{1Prp(LS)}
\mathrm{Ent}_{\sL_{\ti Y_t}}(f^2)\leq C_{LS,\ti Y_t}\int_\R|f'|^2\d\sL_{\ti Y_t}
\end{align}
holds for all $\sL_{\ti Y_t}$-integrable and differentiable function $f:\R\ra\R$,
where
\begin{align*}
C_{LS,\ti Y_t}=2V_T(L_g+L_f(V_T-t))^2\e^{2L_f(V_T-t)}.
\end{align*}
\end{prp}

Here the entropy of $0\leq F\in L^1(\mu)$ with respect to the probability measure $\mu$ is defined as
\begin{align*}
\mathrm{Ent}_{\mu}(F)=\int_\R F\log F\d\mu-\int_\R F\d\mu\cdot\log\int_\R F\d\mu.
\end{align*}
When $\mu$ satisfies \eqref{1Prp(LS)} for $\mu$ replacing $\sL_{\ti Y_t}$, we shall say that $\mu$ satisfies the $LSI(C_{LS})$.
This inequality, initiated by Gross \cite{Gross75}, has become a crucial tool in infinite dimensional stochastic analysis.
It had been well investigated in the context of forward diffusions, and was related with the 2-transportation inequality (see for instance \cite{CGW10,GW06,Ledoux99,OV00,Wang01,Wang09}).

To prove Proposition \ref{Prp(LS)}, we recall a result which asserts that the Logarithmic-Sobolev inequality satisfies stability under push-forward by Lipschitz maps (see \cite[Section 1]{CFJ09} or \cite[Lemma 6.1]{BT20}).

\beg{lem}\label{LS-Le}
If $\psi:(\ti\Om,\|\cdot\|_\infty)\ra\R^d$ is Lipschitzian, i.e.
\begin{align*}
|\psi(\om_1)-\psi(\om_2)|\leq L_\psi\|\om_1-\om_2\|_\infty:=L_\psi\sup\limits_{r\in[0,V_T]}|\om_1(r)-\om_2(r)|,\ \ \om_1, \om_2\in\ti\Om,
\end{align*}
where $L_\psi>0$ is a constant.
Then $\ti\P\circ\psi^{-1}$ satisfies LSI($2V_TL_\psi^2$).
\end{lem}

In light of the proof of Proposition \ref{Prp(Tr)} (see the Appendix), we have shown that
$\ti Y:\ti\Om\ra C([0,V_T])$ is Lipschitz continuous with Lipschitzian constant $(L_g+L_fV_T)\e^{L_fV_T}$,
which actually implies that for any $t\in[0,V_T], \ti Y_t:\ti\Om\ra\R$  is also Lipschitz continuous with Lipschitzian constant $(L_g+L_f(V_T-t))\e^{L_f(V_T-t)}$ thanks to \eqref{add4PfPrp(Tr)} .
So, owing to Lemma \ref{LS-Le}, we get the desired assertion stated in Proposition \ref{Prp(LS)}.

Observe that by Theorem \ref{Th1}, the law of $Y_t$ is the same as that of  $\ti Y_{V_t}$ for every $t\in[0,T]$.
Therefore, by Proposition \ref{Prp(LS)} we have the following Logarithmic-Sobolev inequality for \eqref{Bsde}.

\beg{thm}\label{Th(LS)}
Assume that \textsc{\textbf{(H1)}} holds. Then for any $t\in[0,T]$, the law of $Y_t$ satisfies the LSI($C_{LS,Y_t}$),
where
\begin{align*}
C_{LS,Y_t}=2V_T(L_g+L_f(V_T-V_t))^2\e^{2L_f(V_T-V_t)}.
\end{align*}
\end{thm}

\section{Appendix: proof of Proposition \ref{Prp(Tr)}}

In order to prove the proposition, we first present three lemmas that are needed later on.
The first one concerns the stability of transportation inequalities under push-forward by Lipschitz maps, which is due to \cite[Lemma 2.1]{DGW04} (see, e.g., \cite[Lemma 4.1]{Rid17} and \cite[Corollary 2.2]{SYZ22} for a generalisation).

\beg{lem}\label{FI-Le2}
Let $(E,d_E)$ and $(\bar{E},d_{\bar{E}})$ be two metric spaces.
Assume that $\mu\in T_p(C)$ on $(E,d_E)$ and $\chi:(E,d_E)\ra(\bar{E},d_{\bar{E}})$ is Lipschitz continuous with Lipschitzian constant $L_\chi$.
Then $\mu\circ\chi^{-1}\in T_p(CL^2_\chi)$ on $(\bar{E},d_{\bar{E}})$.
\end{lem}

Our second lemma below provides a sufficient condition expressed in terms of exponential moment for a probability measure satisfying transportation
inequality of the form \eqref{1FI-Le3}.

\beg{lem}\label{FI-Le3}
(\cite[Corollary 2.4]{BV05}) Let $(E,d_E)$ be a metric space, and let $\nu$ be a probability measure on $E$ and $p\geq 1$.
Assume that there exist $x_0\in E$ and $\al>0$ such that
\begin{align*}
\int_E\exp\left\{\al d^{2p}_E(x_0,x)\right\}\d\nu(x)<\infty.
\end{align*}
Then
\begin{align}\label{1FI-Le3}
\W_p(\mu,\nu)\leq C(H(\mu|\nu))^{\ff 1 {2p}},\ \ \mu\in\sP(E)
\end{align}
holds with
\begin{align*}
C=2\inf_{x_0\in E,\al>0}\left[\ff 1{2\al}\left(1+\log\int_E\exp\left\{\al d^{2p}_E(x_0,x)\right\}\d\nu(x)\right)\right]^{\ff 1{2p}}<\infty.
\end{align*}
\end{lem}

Before stating the third lemma, we need some notations from \cite{BT20,EKTZ14}.
For $t\in[0,V_T]$, let $\ti\Om^t$ be the shifted space of $\ti\Om$ given by
\begin{align*}
\ti\Om^t:=\{\ga\in C([t,V_T]):\ga(t)=0\}.
\end{align*}
Denote by $\ti W^t$ and $\ti\P^t$ the canonical process and the Wiener measure on $\ti\Om^t$, respectively,
and by $(\ti\sF^t_s)_{s\in[t,V_T]}$ the filtration generated by  $\ti W^t$.
For $\om\in\ti\Om,t\in[0,V_T]$ and $\ga\in\ti\Om^t$, define the concatenation $\om\otimes_t\ga\in\ti\Om$ by
\begin{equation*}
(\om\otimes_t\ga)(s):=\left\{
\begin{array}{ll}
\om(s),\ \ \ \ \ \ \ \ \ \ \ s\in[0,t),\\
\om(t)+\ga(s),\ \ s\in[t,V_T],
\end{array} \right.
\end{equation*}
and for $\zeta:\ti\Om\times[0,V_T]\ra\R$, define its shift $X^{t,\om}$ as follows
\begin{align*}
\zeta^{t,\om}:&\ \ti\Om^t\times[t,V_T]\ra\R,\cr
&(\ga,s)\mapsto \zeta_s(\om\otimes_t\ga)=:\zeta^{t,\om}_s(\ga).
\end{align*}
As pointed out in \cite{BT20}, we have
\begin{align*}
\ti\E(\zeta|\ti\sF_t)(\om)=\int_{\Om^t}\zeta^{t,\om}(\ga)\ti\P^t(\d\ga)=:\ti\E_{\ti\P^t}\zeta^{t,\om}.
\end{align*}

The following lemma gives a result which is a distribution dependent version of \cite[Lemma 2.2]{BT20}. The proof is pretty similar to that of \cite[Lemma 2.2]{BT20} and we omit it here.

\beg{lem}\label{FI-Le1}
Let $(\ti{Y},\ti{Z})\in\tilde{\Upsilon}$ be the solution to \eqref{Aueq}.
Then for any $t\in[0,V_T]$, there exists a $\ti\P$-zero set $N\subset\ti\Omega$ such that for $\om\in N^c$,
\begin{align*}
\ti{Y}_s^{t,\om}=&g(\ti{W}^{t,\om}_{V_T},\sL_{\ti{W}_{V_T}})
+\int_s^{V_T}f(U_r,\ti{W}^{t,\om}_r,\ti{Y}^{t,\om}_r,\ti{Z}^{t,\om}_r,\sL_{(\ti{W}_r,\tilde{Y}_r,\ti{Y}_r)})\d r\cr
&-\int_s^{V_T}\ti{Z}^{t,\om}_r\d\ti{W}^t_r,
\ \ \ti\P^t \textit{-}a.s., \ s\in[t,V_T]
\end{align*}
and $\ti{Y}_t^{t,\om}=\ti{Y}_t(\om), \ti\P^t$-a.s..
\end{lem}

We are now ready to prove Proposition \ref{Prp(Tr)}.

\emph{Proof of Proposition \ref{Prp(Tr)}.}
Owing to  \textsc{\textbf{(H1)}}, we know that there exists a unique solution  $(\ti{Y},\ti{Z})\in\ti{\Upsilon}$ to \eqref{Aueq}.
Moreover, it easily follows that $\ti{Y}$ has $\ti\P$-almost surely continuous paths and $\ti{Z}$ is square integrable.
The rest of the proof is divided into two steps.

\textsl{Step 1. Transportation inequality for $\ti Y$}.
Using arguments from the proofs of \cite[Proposition 5.4]{EKTZ14} and \cite[Theorem 1.3]{BT20}, we intend to show that $\ti Y:\ti\Om\ra\R$ is Lipschitz continuous.

Let $t\in[0,V_T]$.
According to Lemma \ref{FI-Le1}, there exists a $\ti\P$-zero set $N\subset\ti\Om$ such that for $\omega\in N^c$,
\begin{align}\label{1PfPrp(Tr)}
\ti{Y}_t^{t,\om}=\ti{Y}_t(\om),\ \  \ti\P^t\textit{-}a.s.
\end{align}
and
\begin{align}\label{2PfPrp(Tr)}
\ti{Y}_s^{t,\om}=&g(\ti{W}^{t,\om}_{V_T},\sL_{\ti{W}_{V_T}})
+\int_s^{V_T}f(U_r,\ti{W}^{t,\om}_r,\ti{Y}^{t,\om}_r,\ti{Z}^{t,\om}_r,\sL_{(\ti{W}_r,\tilde{Y}_r,\ti{Y}_r)})\d r\cr
&-\int_s^{V_T}\ti{Z}^{t,\om}_r\d\ti{W}^t_r,
\ \ \ti\P^t\textit{-}a.s., \  s\in[t,V_T].
\end{align}
Then for any $w_1,w_2\in N^c$, by \eqref{2PfPrp(Tr)} and \textsc{\textbf{(H1)}} we derive
\begin{align*}
&\ti{Y}_s^{t,\om_1}-\ti{Y}_s^{t,\om_2}\cr
=&g(\ti{W}^{t,\om_1}_{V_T},\sL_{\ti{W}_{V_T}})-g(\ti{W}^{t,\om_2}_{V_T},\sL_{\ti{W}_{V_T}})\cr
&+\int_s^{V_T}\left(f(U_r,\ti{W}^{t,\om_1}_r,\ti{Y}^{t,\om_1}_r,\ti{Z}^{t,\om_1}_r,\sL_{(\ti{W}_r,\tilde{Y}_r,\ti{Y}_r)})
-f(U_r,\ti{W}^{t,\om_2}_r,\ti{Y}^{t,\om_2}_r,\ti{Z}^{t,\om_2}_r,\sL_{(\ti{W}_r,\tilde{Y}_r,\ti{Y}_r)})\right)
\d r\cr
&-\int_s^{V_T}\left(\ti{Z}^{t,\om_1}_r-\ti{Z}^{t,\om_2}_r\right)\d\ti{W}^t_r\cr
=&g(\ti{W}^{t,\om_1}_{V_T},\sL_{\ti{W}_{V_T}})-g(\ti{W}^{t,\om_2}_{V_T},\sL_{\ti{W}_{V_T}})\cr
&+\int_s^{V_T}\left[\al_r(\ti{W}^{t,\om_1}_r-\ti{W}^{t,\om_2}_r)+\be_r(\ti{Y}^{t,\om_1}_r-\ti{Y}^{t,\om_2}_r)+\rho_r(\ti{Z}^{t,\om_1}_r-\ti{Z}^{t,\om_2}_r) \right]\d r\cr
&-\int_s^{V_T}\left(\ti{Z}^{t,\om_1}_r-\ti{Z}^{t,\om_2}_r\right)\d\ti{W}^t_r,\ \ \P^t\textit{-}a.s., \  s\in[t,V_T],
\end{align*}
where
\begin{align*}
\al_r&:=\int_0^1\partial_xf(U_r,\ti{W}^{t,\om_2}_r+\th(\ti{W}^{t,\om_1}_r-\ti{W}^{t,\om_2}_r),\ti{Y}^{t,\om_1}_r,\ti{Z}^{t,\om_1}_r,\sL_{(\ti{W}_r,\tilde{Y}_r,\ti{Y}_r)})\d\th,\cr
\be_r&:=\int_0^1\partial_yf(U_r,\ti{W}^{t,\om_2}_r,\ti{Y}^{t,\om_2}_r+\th(\ti{Y}^{t,\om_1}_r-\ti{Y}^{t,\om_2}_r),\ti{Z}^{t,\om_1}_r,\sL_{(\ti{W}_r,\tilde{Y}_r,\ti{Y}_r)})\d\th,\cr
\rho_r&:=\int_0^1\partial_zf(U_r,\ti{W}^{t,\om_2}_r,\ti{Y}^{t,\om_2}_r,\ti{Z}^{t,\om_2}_r+\th(\ti{Z}^{t,\om_1}_r-\ti{Z}^{t,\om_2}_r),\sL_{(\ti{W}_r,\tilde{Y}_r,\ti{Y}_r)})\d\th.
\end{align*}
We claim that the product $\e^{\int_t^s\be_r\d r}(\ti{Y}_s^{t,\om_1}-\ti{Y}_s^{t,\om_2})$ yields a more suitable representation for $\ti{Y}_s^{t,\om_1}-\ti{Y}_s^{t,\om_2}$.
Indeed, for $s\in(t,V_T]$,
\begin{align*}
&\d\left[\e^{\int_t^s\be_r\d r}(\ti{Y}_s^{t,\om_1}-\ti{Y}_s^{t,\om_2})\right]\cr
=&\bigg[(\ti{Y}_s^{t,\om_1}-\ti{Y}_s^{t,\om_2})\e^{\int_t^s\be_r\d r}\be_s\cr
&-\e^{\int_t^s\be_r\d r}
\left(\al_s(\ti{W}^{t,\om_1}_s-\ti{W}^{t,\om_2}_s)+\be_s(\ti{Y}^{t,\om_1}_s-\ti{Y}^{t,\om_2}_s)+\rho_s(\ti{Z}^{t,\om_1}_s-\ti{Z}^{t,\om_2}_s)\right)\bigg]\d s\cr
&+e^{\int_t^s\be_r\d r}\left(\ti{Z}^{t,\om_1}_s-\ti{Z}^{t,\om_2}_s\right)\d\ti{W}^t_s\cr
=&\left[
-\e^{\int_t^s\be_r\d r}
\left(\al_s(\ti{W}^{t,\om_1}_s-\ti{W}^{t,\om_2}_s)+\rho_s(\ti{Z}^{t,\om_1}_s-\ti{Z}^{t,\om_2}_s)\right)\right]\d s
+e^{\int_t^s\be_r\d r}\left(\ti{Z}^{t,\om_1}_s-\ti{Z}^{t,\om_2}_s\right)\d\ti{W}^t_s.
\end{align*}
Then, integrating from $s$ to $V_T$ we get
\begin{align*}
&\e^{\int_t^{V_T}\be_r\d r}(\ti{Y}_{V_T}^{t,\om_1}-\ti{Y}_{V_T}^{t,\om_2})-\e^{\int_t^s\be_r\d r}(\ti{Y}_s^{t,\om_1}-\ti{Y}_s^{t,\om_2})\cr
=&-\int_s^{V_T}\e^{\int_t^r\be_\th\d \th}\left(\al_r(\ti{W}^{t,\om_1}_r-\ti{W}^{t,\om_2}_r)+\rho_r(\ti{Z}^{t,\om_1}_r-\ti{Z}^{t,\om_2}_r)\right)\d r\cr
&+\int_s^{V_T}e^{\int_t^r\be_\th\d\th}\left(\ti{Z}^{t,\om_1}_r-\ti{Z}^{t,\om_2}_r\right)\d\ti{W}^t_r,\ \ \ti\P^t\textit{-}a.s..
\end{align*}
Observing that $\ti{Y}_{V_T}^{t,\om_i}=g(\ti{W}^{t,\om_i}_{V_T},\sL_{\ti{W}_{V_T}}),i=1,2$, we thus obtain
\begin{align}\label{3PfPrp(Tr)}
\ti{Y}_s^{t,\om_1}-\ti{Y}_s^{t,\om_2}
=&\e^{\int_s^{V_T}\be_r\d r}\left(g(\ti{W}^{t,\om_1}_{V_T},\sL_{\ti{W}_{V_T}})-g(\ti{W}^{t,\om_2}_{V_T},\sL_{\ti{W}_{V_T}})\right)\cr
&+\int_s^{V_T}\e^{\int_s^r\be_\th\d \th}\left(\al_r(\ti{W}^{t,\om_1}_r-\ti{W}^{t,\om_2}_r)+\rho_r(\ti{Z}^{t,\om_1}_r-\ti{Z}^{t,\om_2}_r)\right)\d r\cr
&-\int_s^{V_T}e^{\int_s^r\be_\th\d\th}\left(\ti{Z}^{t,\om_1}_r-\ti{Z}^{t,\om_2}_r\right)\d\ti{W}^t_r\cr
=&\e^{\int_s^{V_T}\be_r\d r}\left(g(\ti{W}^{t,\om_1}_{V_T},\sL_{\ti{W}_{V_T}})-g(\ti{W}^{t,\om_2}_{V_T},\sL_{\ti{W}_{V_T}})\right)\nonumber\\
&+\int_s^{V_T}\e^{\int_s^r\be_\th\d \th}\al_r(\ti{W}^{t,\om_1}_r-\ti{W}^{t,\om_2}_r)\d r
-\int_s^{V_T}e^{\int_s^r\be_\th\d\th}\left(\ti{Z}^{t,\om_1}_r-\ti{Z}^{t,\om_2}_r\right)\d\bar{W}^t_r,
\end{align}
where $\bar{W}^t_r:=\ti W_r-\int_t^r\rho_r\d r$.\\
Now, for $s\in[t,V_T]$, we set
\begin{align*}
R_s=\exp\left\{\int_t^s\rho_r\d\ti W_r-\ff 1 2 \int_t^s|\rho_r|^2\d r\right\}.
\end{align*}
By \textsc{\textbf{(H1)}}, it is easy to verify that the Novikov condition holds, which implies that $\bar{W}^t_\cdot$ is a Brownian motion under the probability $R_{V_T}\ti\P^t$ due to the Girsanov theorem.
Then, conditioning by $\ti\sF^t_s$ under $R_{V_T}\ti\P^t$ on both sides of \eqref{3PfPrp(Tr)} yields
\begin{align}\label{4PfPrp(Tr)}
&\ti\E_{R_{V_T}\ti\P^t}\left(\ti{Y}_s^{t,\om_1}-\ti{Y}_s^{t,\om_2}|\ti\sF^t_s\right)\cr
=&\ti\E_{R_{V_T}\ti\P^t}\left[\e^{\int_s^{V_T}\be_r\d r}\left(g(\ti{W}^{t,\om_1}_{V_T},\sL_{\ti{W}_{V_T}})-g(\ti{W}^{t,\om_2}_{V_T},\sL_{\ti{W}_{V_T}})\right)|\ti\sF^t_s\right]\nonumber\\
&+\ti\E_{R_{V_T}\ti\P^t}\left[\int_s^{V_T}\e^{\int_s^r\be_\th\d \th}\al_r(\ti{W}^{t,\om_1}_r-\ti{W}^{t,\om_2}_r)\d r|\ti\sF^t_s\right],\ \ \ti\P^t\textit{-}a.s., \  s\in[t,V_T].
\end{align}
Consequently, we have
\begin{align*}
\ti{Y}_t^{t,\om_1}-\ti{Y}_t^{t,\om_2}
=&\ti\E_{R_{V_T}\ti\P^t}\left[\e^{\int_t^{V_T}\be_r\d r}\left(g(\ti{W}^{t,\om_1}_{V_T},\sL_{\ti{W}_{V_T}})-g(\ti{W}^{t,\om_2}_{V_T},\sL_{\ti{W}_{V_T}})\right)\right]\cr
&+\ti\E_{R_{V_T}\ti\P^t}\left[\int_t^{V_T}\e^{\int_t^r\be_\th\d \th}\al_r(\ti{W}^{t,\om_1}_r-\ti{W}^{t,\om_2}_r)\d r\right],\ \ \ti\P^t\textit{-}a.s..
\end{align*}
This allows us to deduce from \eqref{1PfPrp(Tr)} and \textsc{\textbf{(H1)}} that
\begin{align}\label{add4PfPrp(Tr)}
&|\ti{Y}_t(\om_1)-\ti{Y}_t(\om_2)|=|\ti{Y}_t^{t,\om_1}-\ti{Y}_t^{t,\om_2}|\cr
\leq&\e^{L_f(V_T-t)}\left(L_g\ti\E_{R_{V_T}\ti\P^t}|\ti{W}^{t,\om_1}_{V_T}-\ti{W}^{t,\om_2}_{V_T}|
+L_f\int_t^{V_T}\ti\E_{R_{V_T}\ti\P^t}|\ti{W}^{t,\om_1}_r-\ti{W}^{t,\om_2}_r|\d r\right)\cr
=&\e^{L_f(V_T-t)}\bigg(L_g\ti\E_{R_{V_T}\ti\P^t}|(\om_1\otimes_t\cdot)(V_T)-(\om_2\otimes_t\cdot)(V_T)|\cr
&\ \ \ \ \ \ \ \ \ \ \ \ \ +L_f\int_t^{V_T}\ti\E_{R_{V_T}\ti\P^t}|(\om_1\otimes_t\cdot)(r)-(\om_2\otimes_t\cdot)(r)|\d r\bigg)\cr
\leq&\left(L_g+L_f(V_T-t)\right)\e^{L_f(V_T-t)}\sup\limits_{0\leq r\leq t}|\om_1(r)-\om_2(r)|,
\end{align}
where the last inequality is due to the definitions of concatenation variables $\om_i\otimes_t\cdot, i=1,2$.
Noting that $t\in[0,V_t]$ and $w_1,w_2\in N^c$ are arbitrary and $\ti\P(N)=0$, we conclude that $\ti Y:\ti\Om\ra C([0,V_T])$ is Lipschitz continuous with Lipschitzian constant $(L_g+L_fV_T)\e^{L_fV_T}$.
Therefore, taking into account of the fact that the law of Wiener process satisfies $T_2(2)$ (see \cite[Theorem 3.1]{FU04}),
we obtain the first assertion due to Lemma \ref{FI-Le2}.

\textsl{Step 2. Transportation inequality for $\ti Z$}.
We first suppose that $g(x,\mu)$ and $f(t,x,y,z,\nu)$ are differentiable with respect to $x,y$ and $z$.
Then $(\ti Y,\ti Z)$  is differentiable, and moreover $(\nabla\ti{Y},\nabla\ti{Z})$ solves the following linear DDBSDE
\begin{align*}
\nabla\ti{Y}_t=&\nabla_xg(\tilde{W}_{V_T},\sL_{\tilde{W}_{V_T}})-\int_t^{V_T}\nabla\tilde{Z}_s\d\tilde{W}_s\cr
&+\int_t^{V_T}\Big[\nabla_xf(U_s,\tilde{W}_s,\tilde{Y}_s,\tilde{Z}_s,\sL_{(\tilde{W}_s,\tilde{Y}_s,\tilde{Z}_s)})
+\nabla_yf(U_s,\tilde{W}_s,\tilde{Y}_s,\tilde{Z}_s,\sL_{(\tilde{W}_s,\tilde{Y}_s,\tilde{Z}_s)})\nabla\ti{Y}_s\cr
&\qquad\qquad+\nabla_zf(U_s,\tilde{W}_s,\tilde{Y}_s,\tilde{Z}_s,\sL_{(\tilde{W}_s,\tilde{Y}_s,\tilde{Z}_s)})\nabla\ti{Z}_s\Big]\d s.
\end{align*}
Along the same lines as in \eqref{4PfPrp(Tr)}, applying the product $\e^{\int_0^t\nabla_yf(U_s,\tilde{W}_s,\tilde{Y}_s,\tilde{Z}_s,\sL_{(\tilde{W}_s,\tilde{Y}_s,\tilde{Z}_s)})\d s}\nabla\ti{Y}_t$
and the Girsanov theorem we deduce that there exists some probability $\ti\Q$ under which
\begin{align*}
\tilde{W}_\cdot-\int_0^\cdot\nabla_zf(U_s,\tilde{W}_s,\tilde{Y}_s,\tilde{Z}_s,\sL_{(\tilde{W}_s,\tilde{Y}_s,\tilde{Z}_s)})\d s
\end{align*}
is a Brownian motion, and $\nabla\ti{Y}_t$ has the representation
\begin{align*}
\nabla\ti{Y}_t=&\ti\E_{\ti\Q}\bigg[e^{\int_t^{V_T}\nabla_yf(U_s,\tilde{W}_s,\tilde{Y}_s,\tilde{Z}_s,\sL_{(\tilde{W}_s,\tilde{Y}_s,\tilde{Z}_s)})\d s}\nabla_xg(\tilde{W}_{V_T},\sL_{\tilde{W}_{V_T}})\cr
&+\int_t^{V_T}e^{\int_t^s\nabla_yf(U_r,\tilde{W}_r,\tilde{Y}_r,\tilde{Z}_r,\sL_{(\tilde{W}_r,\tilde{Y}_r,\tilde{Z}_r)})\d r}\nabla_xf(U_s,\tilde{W}_s,\tilde{Y}_s,\tilde{Z}_s,\sL_{(\tilde{W}_s,\tilde{Y}_s,\tilde{Z}_s)})\d s\big|\ti\sF_t\bigg].
\end{align*}
Consequently, it is easy to see that
\begin{align}\label{5PfPrp(Tr)}
|\nabla\ti{Y}_t|\leq\e^{L_f(V_T-t)}(L_g+L_f(V_T-t)).
\end{align}
On the other hand, it is classical to show that there exists a version of $(\ti Z_t)_{t\in[0,T]}$ given by $(\nabla\ti{Y}_t)_{t\in[0,T]}$
(see, e.g., \cite[Remark 9.1]{Li18}).
Hence, combining this with \eqref{5PfPrp(Tr)} leads to
\begin{align}\label{6PfPrp(Tr)}
|\ti Z_t|\leq\e^{L_f(V_T-t)}(L_g+L_f(V_T-t)).
\end{align}
When $g$ and $f$ are not differentiable, we can also obtain \eqref{6PfPrp(Tr)} via a standard approximation and stability results (see Remark \ref{Re1} (i)).
Therefore, with the help of Lemma \ref{FI-Le3}, we derive the second assertion.
\qed

\textbf{Acknowledgement}

X. Fan is partially supported by the Natural Science Foundation of Anhui Province (No. 2008085MA10) and the National Natural Science Foundation of China (No. 11871076, 12071003).

\end{document}